\documentclass[11pt]{article}

\usepackage[margin=1in]{geometry}
\usepackage{authblk}
\usepackage{amsmath,amssymb,amsfonts,bm}
\usepackage{graphicx}
\graphicspath{{./}{advection_aligned_eval_figs/}{advection_ablation_figs/}{advection_frequency8_run/figures/}{burgers_train_true_spectrum_figs/}{burgers_eval_figs/}{burgers_corrupted_context_query_run/wrong_eval_figs/}{burgers_corrupted_context_query_run/figs/}{burgers_wrong_eval_figs/}{allen_cahn_eval_figs/}{allen_cahn_aligned_eval_figs/}{allen_cahn_frequency8_run/figures/}}
\usepackage{subcaption}
\usepackage{booktabs}
\usepackage{tabularx}
\usepackage{array}
\usepackage{makecell}
\usepackage{microtype}
\usepackage{enumitem}
\usepackage{xcolor}
\usepackage[numbers,sort&compress]{natbib}
\usepackage{hyperref}
\usepackage{placeins}
\usepackage{float}
\usepackage[font=small,labelfont=bf]{caption}

\newcolumntype{Y}{>{\raggedright\arraybackslash}X}
\newcolumntype{C}{>{\centering\arraybackslash}X}

\hypersetup{
    colorlinks=true,
    linkcolor=blue,
    citecolor=blue,
    urlcolor=blue
}

\newcommand{\cC}{\mathcal{C}}
\newcommand{\cA}{\mathcal{A}}
\newcommand{\cU}{\mathcal{U}}
\newcommand{\cP}{\mathcal{P}}
\newcommand{\uq}{u_{\mathrm{q}}}

\newcommand{\GICON}{G_{\cC}^{\mathrm{ICON}}}

\newcommand{\Err}{\mathrm{Err}}
\newcommand{\diag}{\operatorname{diag}}

\newcommand{\JICON}{J_{\cC}^{\mathrm{ICON}}}
\newcommand{\Jtrue}{J^{\dagger}}
\newcommand{\AICON}{A_{\cC}^{\mathrm{ICON}}}
\newcommand{\Atrue}{A^{\dagger}}

\newcommand{\dd}{\mathrm{d}}

\title{Spectral Audit of In-Context Operator Networks}
\author[1]{Zhiwei Gao}
\author[2]{Liu Yang}
\author[1,*]{George Em Karniadakis}

\affil[1]{Division of Applied Mathematics, Brown University, Providence, RI, USA}
\affil[2]{Department of Mathematics, National University of Singapore, Singapore}
\affil[*]{Corresponding author: \texttt{george\_karniadakis@brown.edu}}

\date{}

\begin{document}
\maketitle

\begin{abstract}

Existing evaluations of neural operators and in-context operator learning rely primarily on prediction error, but accurate output prediction does not guarantee the correct local dynamical structure. A model may match solutions while exhibiting incorrect sensitivities, distorted frequency response, spurious mode coupling, or unstable tangent behavior.
We introduce a Jacobian-based spectral audit for in-context operator learning. For a fixed prompt, we differentiate the network output with respect to the query function and view the resulting Jacobian as a learned tangent operator. Projecting it onto Fourier modes, we obtain a local spectral characterization of the inferred operator, including frequency-dependent gains, phase structure, and cross-mode coupling.
The audit complements standard prediction metrics by testing whether the model reproduces local mechanisms of the underlying PDE operator rather than only outputs. Across benchmarks, the audit reveals distinct operator-level phenomena, including phase transport, viscosity-dependent damping, nonlinear mode coupling, and reaction--diffusion stability structure. It also detects failures partially hidden by prediction-error metrics, including high-frequency degradation, incorrect phase recovery, and prompt--operator inconsistencies. Corrupted or internally inconsistent prompts lead to degraded tangent-operator structure even when pointwise predictions remain partially accurate.
Our results suggest that prediction accuracy and local operator fidelity are distinct properties of learned neural operators. Our framework also provides a diagnostic for stability, sensitivity, and operator consistency.

\end{abstract}

\iffalse
\begin{abstract}
In-context operator networks infer an operator from a small prompt and apply it to new query functions. Standard evaluation mainly relies on prediction error, but accurate output prediction does not necessarily imply that the inferred operator has the correct local structure. We propose a Jacobian-based spectral audit for in-context operator learning. For a fixed prompt, we view the trained network as a local map from query inputs to predicted outputs, differentiate this map with respect to the query function, and analyze the resulting tangent operator in a Fourier basis. This diagnostic complements pointwise prediction metrics by revealing structural errors that may otherwise remain hidden, including loss of high-frequency response, spurious mode coupling, and collapse toward an averaged operator. The proposed audit provides a practical way to assess whether an in-context model has learned a meaningful operator representation and to compare the effects of training data, prompt quality, and model design beyond prediction error alone.
\end{abstract}
\fi 

% ------------------------------------------------------------

\section{Introduction}

Neural operators provide a framework for learning maps between function spaces, such as solution operators of parametrized partial differential equations (PDE)~\cite{lu2021learning,kovachki2023neural}. Representative architectures, including DeepONet, Fourier Neural Operator, transformer-based neural operators, and other variants \cite{cao2024laplace, tripura2022wavelet}, have become standard tools for approximating PDE solution maps from data~\cite{mao2021deepm,nath2026digital, bodnar2025foundation}. In many applications, however, the operator to be applied at test time is not specified directly by an explicit parameter. Instead, the model is given a small number of input--output examples from the current physical regime and must use this prompt to predict the response to a new query. This is the setting of {\em in-context operator learning}. In the case of In-Context Operator Networks (ICONs), the hidden PDE parameter is not supplied as an input; the network must infer the relevant operator from the context pairs~\cite{yang2023context,yang2024pde}.

The usual way to evaluate such a model is to measure the prediction error on query outputs. This is necessary, but it does not fully determine whether the model has learned the correct operator. A context-conditioned model may produce accurate predictions on a finite set of query functions while responding incorrectly to small perturbations of the query. For PDE solution operators, this local response is meaningful: it describes the tangent map of the dynamics and encodes stability, frequency-dependent damping or amplification, and interactions between modes. These properties become important when learned operators are used beyond one-step prediction, for example in rollouts, inverse problems, uncertainty propagation, or design and optimization.

This viewpoint is related to a broader theme in scientific machine learning
(SciML): matching pointwise data is often insufficient when the learned model is
intended to represent a physical system. Physics-informed methods incorporate
differential-equation structure into training objectives
~\cite{raissi2019physics,sirignano2018dgm,karniadakis2021physics}, while
data-driven equation discovery aims to identify governing equations or active
terms from observations~\cite{brunton2016discovering,rudy2017data}. Our goal is
different. We do not propose a new neural-operator architecture, add a PDE
residual to the training loss, or attempt to discover an equation from data.
Instead, we take a trained context-conditioned neural operator and examine the
local operator it has inferred once a prompt is fixed. This also differs from
imposing spectral constraints during training~\cite{behroozi2025sensitivity}
or analyzing the sensitivity of a single fixed neural operator
~\cite{shikhman2026forcing}. Rather, we provide a posterior spectral analysis
of the inferred in-context neural operator in a continuous, operator-theoretic
sense.

For a given context $\cC$, ICON induces a deterministic map from the query function to the predicted output, denoted by $\GICON$. We study the Jacobian $D\GICON(\uq)$ with respect to the query input $\uq$ and interpret it as a learned tangent operator. The main question is whether this learned tangent map agrees with the tangent map of the reference PDE solution operator, not only whether the predicted output is accurate for one query.

To answer this question, we propose a Jacobian-based spectral audit. The Jacobian is computed or probed through directional derivatives and then analyzed in a real Fourier basis. This gives a simple way to inspect how the learned local map responds to perturbations at different frequencies and how it transfers information between modes. We use Fourier-mode gains, accumulated spectral errors, and Fourier-projected feature maps as complementary diagnostics. The Fourier basis is not essential to the method, but it is especially transparent for the one-dimensional periodic examples considered here. It provides a practical protocol for checking whether a trained in-context operator model has learned a meaningful local operator representation. The audit can reveal structural errors that are not visible from average prediction error alone, such as loss of high-frequency response, incorrect mode coupling, distorted phase structure, or collapse toward an averaged operator. Our experiments show that prediction accuracy and local operator consistency are related but distinct, suggesting that Jacobian-based diagnostics are a useful complement to standard error metrics when sensitivity and stability matter.

% ------------------------------------------------------------
\section{In-Context Operator Learning}

Let \(\cA\) and \(\cU\) be spaces of input and output functions.  A standard
neural operator approximates a fixed map
\[
    G^\dagger:\cA\to\cU,
    \qquad u\mapsto G^\dagger(u),
\]
where $G^{\dagger}$ is the single true solution operator. The problems studied here involve a family of solution operators
\[
    \{G_\eta^\dagger:\cA\to\cU\}_{\eta\in\cP},
\]
where \(\eta\) is a hidden physical parameter and $\mathcal{P}$ is the parameter space.  A parameter-conditioned neural operator would
usually receive \(\eta\) as an explicit input. While for ICON, it is given a
small prompt generated by one unknown operator instance and must infer that
operator from the prompt. Formally, the ICON is defined to approximate the map:
\begin{equation}
    \label{ICON}
    F^{\dagger}: \mathfrak{C}\times \mathcal{A}\rightarrow \mathcal{U},
\end{equation}
where $\mathfrak{C}$ is the space of the context. For a fixed \(\eta\), the context is $\cC_\eta=\{(u_i,v_i)\}_{i=1}^{N_{\mathrm{demo}}}, v_i=G_\eta^\dagger(u_i)$ and $N_{\mathrm{demo}}$ is the number of solution pairs we use in the context. 
Given \(\cC_\eta\) and a new query function \(\uq\), the network predicts
\[
    \widehat v_{\mathrm q}=F_\theta(\cC_\eta,\uq)
    \approx F^{\dagger}(\mathcal{C}_{\eta}, u_{q}) =  G_\eta^\dagger(\uq),
\]
where $F_{\theta}$ is the ICON taking input with the context and the query to approximate the target quantity of interest (QoI) $G^{\dagger}(u_{q})$.
Thus, the same weights \(\theta\) are shared across many hidden operators, while
the context selects the effective operator for the current query.

In practice, functions are represented as key--value tokens on the spatial grid.  A condition function contributes tokens
\(\{(x_j,u(x_j))\}_j\), and the corresponding QoI contributes
\(\{(x_j,G_\eta^\dagger u(x_j))\}_j\).  The training objective is
the usual squared error on the output tokens, including the query
QoI tokens and the context QoI tokens used in the ICON construction. 

After fixing a context, the trained model becomes a concrete operator
\[
    G_{\cC_\eta}^{\mathrm{ICON}}:\cA\to\cU,
    \qquad
    G_{\cC_\eta}^{\mathrm{ICON}}(\uq)=F_\theta(\cC_\eta,\uq).
\]
This fixed-context map is the object we audit.  During the audit, the context is
held fixed and only the query is perturbed:
\begin{equation}
\label{ICON_operator}
    \JICON(\uq)=D_{u_{q}} G_{\cC_\eta}^{\mathrm{ICON}}(u_{q}),
\end{equation}
where $D_{u_{q}}$ is the differential operator with respect to the input.
The reference object is the PDE tangent map
\begin{equation}
\label{true_operator}
\Jtrue_\eta(\uq)=D_{u_q}G_\eta^\dagger(\uq),
\end{equation}
which measures how a perturbation of the initial condition changes the final
state under the same physical parameter.  Comparing \(\JICON\) and \(\Jtrue\), 
therefore, tests the local mechanism of the inferred operator, not only its
pointwise output accuracy.

% ------------------------------------------------------------
\section{Experimental Setup and Unified Metrics}
\label{sec:experimental_setup}

In this section, we present the experimental setup for the following tasks. All main experiments use the same learning problem: given a prompt generated by
one hidden PDE parameter, learn the finite-time solution operator
\[
    G_\eta^\dagger:u_0(x)\mapsto u_{T}(x;\eta), x\in [0, 1],
\]
where $u_{T}$ is the solution at time $T$ and $\eta$ is the physical parameter. For different PDEs, $\eta$ has different meanings. For example, for transport equation, $\eta$ represents the velocity. For Burgers' equation, it represents the viscosity. 

\subsection{Data generation}
\label{subsec:data_generation_protocol}
 Each prompt in $C_{\eta}$ contains five solution pairs , i.e., $N_{demo} = 5$,  and one query from the same
hidden operator.  All functions are sampled on the periodic grid
\(x_j=j/n_{\mathrm{grid}}\) with \(n_{\mathrm{grid}}=100\).  Random initial
conditions are drawn from
\[
    u_0(x)=b+
    \sum_{m=1}^{n_{\mathrm{train}}}
    \left(a_m^{\sin}\sin(2\pi m x)+a_m^{\cos}\cos(2\pi m x)\right),
\]
where
\[
    a_m^{\sin},a_m^{\cos}\sim
    \mathcal N\!\left(0,\frac{\sigma_u^2}{m^{2\alpha_u}}\right),
    \qquad b\sim\mathcal N(0,\sigma_b^2).
\]
In practice, we use \(n_{\mathrm{train}}=16\), \(\sigma_u=0.7\),
\(\alpha_u=1\), and \(\sigma_b=0\).   The network is a Transformer ICON based on masked self-attention~\cite{vaswani2017attention}.  The attention mask implements the prompt structure: during training the
model reconstructs QoI values from the visible condition and context information, while during evaluation it outputs the query QoI.  We
use the same architecture and optimizer for the three main PDE examples. The detailed hyperparameter settings are in Table \ref{tab:common_setup_unified}. The details of the true operator are listed in Table \ref{tab:pde_setup_unified}.

\begin{table}[H]
\centering
\small
\caption{Experimental setup for all of the experiments.}
\label{tab:common_setup_unified}
\begin{tabularx}{0.92\textwidth}{lY}
\toprule
Setting & Value \\
\midrule
Prompt & 5 demonstrations + 1 query, 100 spatial tokens per function \\
Architecture & 6 layers, 8 heads, model dimension 256, head dimension 32 \\
Optimization & AdamW~\cite{loshchilov2017decoupled}, batch size 4, 200000 iterations \\
Offline pool & 2000 hidden operators $\times$ 20 input--output pairs/operator \\
\bottomrule
\end{tabularx}
\end{table}

\begin{table}[H]
\centering
\small
\caption{PDE-specific finite-time operators and data-generation settings.}
\label{tab:pde_setup_unified}
\begin{tabularx}{\textwidth}{lYll}
\toprule
Example & Target operator & Hidden parameter & Horizon \\
\midrule
Advection & $G_c^{T}u_0(x)=u_0((x-c T)\bmod 1)$ & $c\in[0.5,1.0]$ & $T=0.5$ \\
Burgers & $G_c^T:u_0\mapsto u(T,\cdot)$ for $u_t+u u_x=c u_{xx}$ & $c\in[0.005,0.05]$ & $T=0.25$ \\
Allen--Cahn & $G_\varepsilon^T:u_0\mapsto u(T,\cdot)$ for $u_t=\varepsilon^2u_{xx}+u-u^3$ & $\varepsilon\in[0.02,0.08]$ & $T=0.25$ \\
\bottomrule
\end{tabularx}
\end{table}

\subsection{Audit quantities and error metrics}
\label{subsec:jacobian_audit_quantities}

For a fixed context \(\cC_\eta\), we compare the learned query Jacobian and the
true PDE tangent map defined in Eqs. \eqref{ICON_operator} and \eqref{true_operator}. To calculate the spectrum, we add Fourier basis functions as the perturbations, 
\[
    \phi_k^{\cos}(x)=\sqrt{2}\cos(2\pi kx),
    \qquad
    \phi_k^{\sin}(x)=\sqrt{2}\sin(2\pi kx).
\]
Then, we can define the spectrum, also called the gains, as 
\[
    \widehat g_\alpha(k;\eta)=
    \frac{\|\JICON(\uq)\phi_k^\alpha\|_2}{\|\phi_k^\alpha\|_2},
    \qquad
    g_\alpha^\dagger(k;\eta)=
    \frac{\|\Jtrue_\eta(\uq)\phi_k^\alpha\|_2}{\|\phi_k^\alpha\|_2},
    \quad \alpha\in\{\cos,\sin\}.
\]
To compare with the ground truth, we defined the accumulated relative error as 
\[
    \Err_\alpha(K;\eta)=
    \frac{\left(\sum_{k=1}^{K}
    |\widehat g_\alpha(k;\eta)-g_\alpha^\dagger(k;\eta)|^2\right)^{1/2}}
    {\left(\sum_{k=1}^{K}|g_\alpha^\dagger(k;\eta)|^2\right)^{1/2}}.
\]
The cumulative form is intentional: in diffusive problems, pointwise relative
errors at high modes can be dominated by tiny denominators.

To visualize mode mixing, we also project the Jacobian onto a finite real
Fourier basis ordered by sine--cosine pairs.  Let
\[
    \Phi
    =
    [\phi^{\cos}_1,\phi^{\sin}_1,\ldots,
    \phi^{\cos}_M,\phi^{\sin}_M],
    \qquad
    \Psi
    =
    [\psi^{\cos}_1,\psi^{\sin}_1,\ldots,
    \psi^{\cos}_M,\psi^{\sin}_M],
\]
where \(\Phi\in\mathbb R^{n_{\rm in}\times 2M}\) and
\(\Psi\in\mathbb R^{n_{\rm out}\times 2M}\) are the input- and output-side
Fourier basis matrices.  The
projected maps are
\[
    \AICON=\Psi^\top\JICON\Phi,
    \qquad
    \Atrue_\eta=\Psi^\top\Jtrue_\eta\Phi.
\]
Thus, each entry records how one input Fourier direction contributes to one
output Fourier direction.  We report the relative Frobenius error
\[
    \Err_A(\eta)=
    \frac{\|\AICON-\Atrue_\eta\|_F}{\|\Atrue_\eta\|_F}.
\]
and use the off-diagonal ratio
\[
    r_{\mathrm{off}}=\frac{\|A-\diag(A)\|_F}{\|A\|_F}
\]
as a compact measure of cross-mode coupling.  Jacobian-vector products and full
Jacobians are obtained by automatic differentiation through the ICON predictor
and through the discrete reference solver.  The continuous variational equations
behind these derivatives are collected in Appendix~\ref{app:true_derivatives}.

\FloatBarrier

% ------------------------------------------------------------

\section{Example I: Hidden-Shift Advection}

The advection experiment is the cleanest starting point because the exact
solution operator is known and has a simple Fourier representation.  We solve
\[
    u_t+c u_x=0,\qquad x\in[0,1],
\]
with periodic boundary conditions and observe the solution after a fixed time
\(T\).  The finite-time operator is defined in Table \ref{tab:pde_setup_unified}. 
Thus the hidden speed \(c\) controls only the translation amount.  In Fourier
space this means that the operator should not smooth, amplify, or transfer
energy between different frequencies.  Its only nontrivial action is to rotate
each sine--cosine pair by a phase determined by \(c T\).

\begin{figure}[!htbp]
    \centering
    \includegraphics[width=0.5\textwidth]{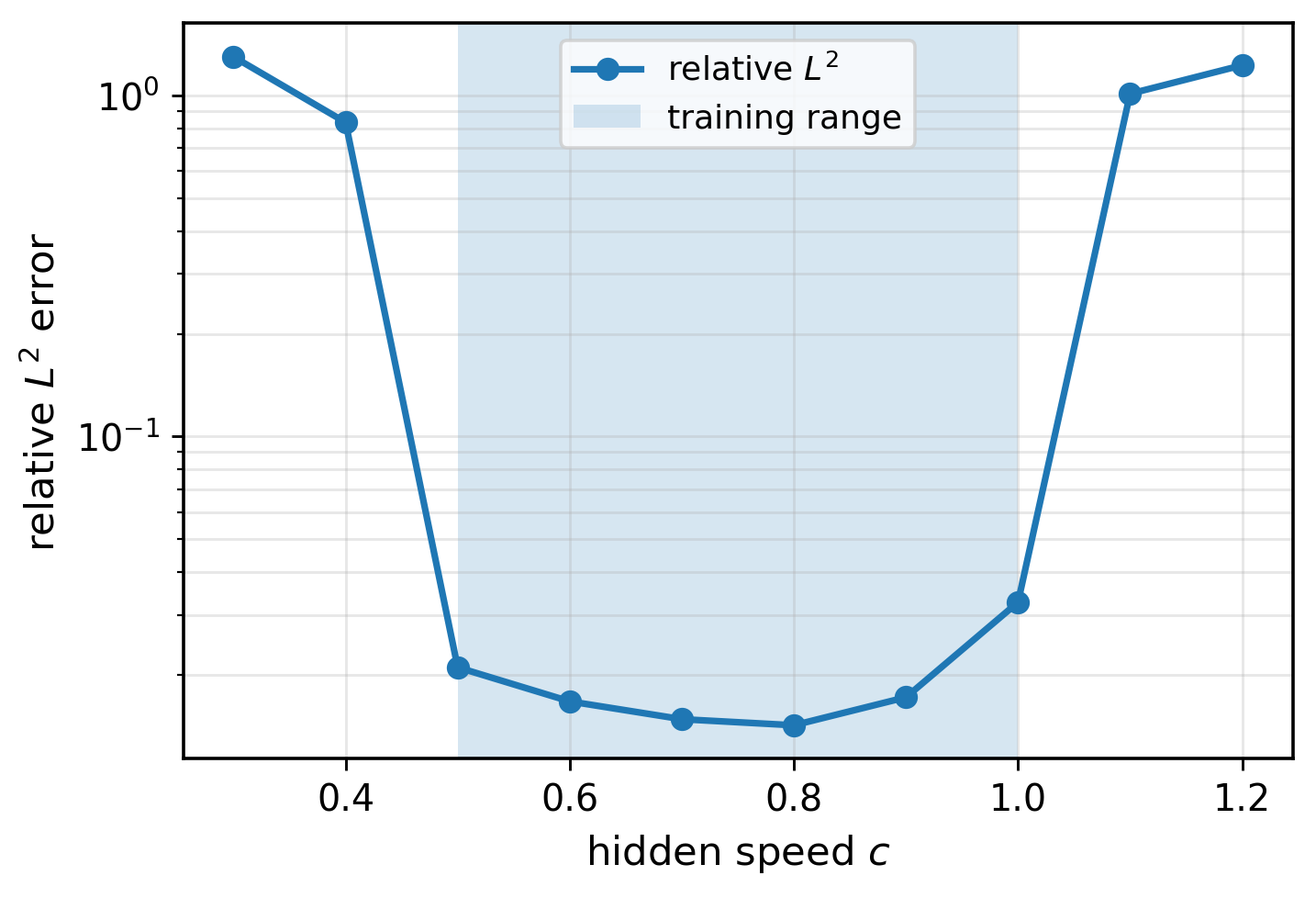}
    \caption{Advection relative \(L^2\) prediction error over hidden speed. The shaded region denotes the training range.}
    \label{fig:adv_rel_l2}
\end{figure}

The prediction-error curve in Figure~\ref{fig:adv_rel_l2} serves as a basic
check that the context is sufficient for ICON to infer the hidden shift.  Within
the training interval, the relative error stays small, while the error increases
as the speed moves away from the trained range.  This is the expected behavior
for a context model trained on a bounded family of shifts: it can interpolate
between seen translation amounts, but extrapolating to larger or smaller shifts
is more difficult.  We therefore interpret the Jacobian diagnostics primarily
inside the shaded interval.

\begin{figure}[!htbp]
    \centering
    \includegraphics[width=0.8\textwidth]{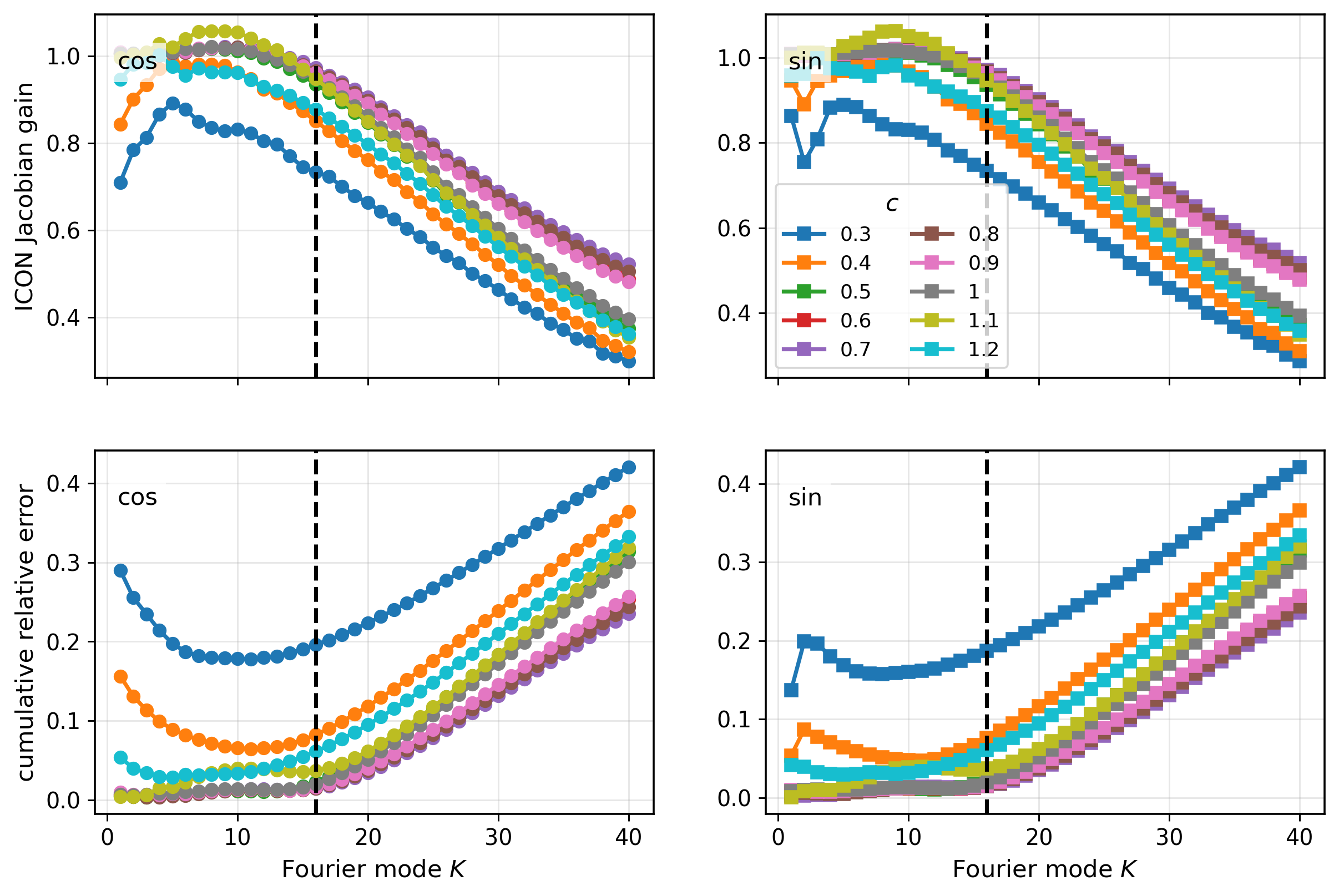}
    \caption{Advection query spectrum over hidden speed. The top row shows ICON Fourier-mode gains, and the bottom row reports accumulated relative spectral errors.}
    \label{fig:adv_query_spectrum}
\end{figure}

Because advection is linear, the query Jacobian should be exactly the same shift
operator for every base state:
\[
    D_{u_0}G_c^{T}(u_0)v(x)=v((x-cT)\bmod1).
\]
This makes the gain diagnostic particularly easy to read.  The reference gain
is 1 for every Fourier mode, so the top row of
Figure~\ref{fig:adv_query_spectrum} asks whether the learned Jacobian preserves
amplitudes.  The curves remaining close to one over the training bandwidth
indicate that the model has not turned transport into an artificial smoothing
operator. While out of the training range, the spectrum decays to zero fast due to the spectral bias of neural networks, which has been studied in many literature \cite{rahaman2019spectral, khodakarami2026spectral}.  The bottom row gives a cumulative relative error.  This quantity is
more stable than a pointwise ratio and shows how the small gain discrepancies
accumulate as more modes are included.

\begin{figure}[!htbp]
    \centering
    \includegraphics[width=0.5\textwidth]{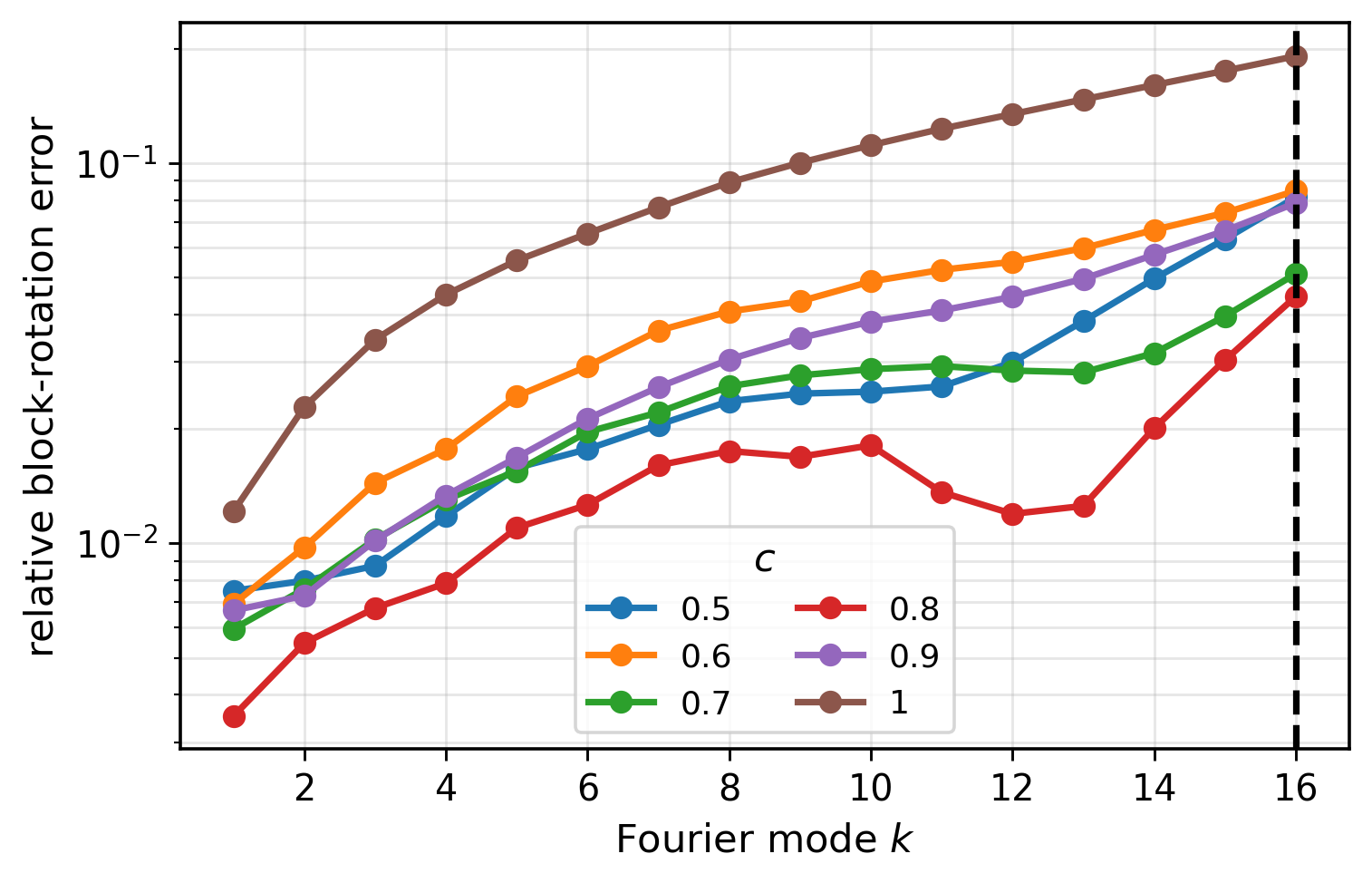}
    \caption{Advection phase and block-rotation error over hidden speed.}
    \label{fig:adv_phase_block_error}
\end{figure}

For transport, however, gain alone is not enough.  A map may preserve the norm
of each Fourier mode while rotating the sine--cosine pair by the wrong angle,
which would correspond to an incorrect shift.  This is why the phase/block error
in Figure~\ref{fig:adv_phase_block_error} is the decisive diagnostic for the
advection case.  In the real Fourier basis, the exact tangent map is a direct
sum of \(2\times2\) rotation blocks.  Small block error therefore means that the
context has identified the translation parameter itself, not just an
amplitude-preserving approximation.

\begin{figure}[!htbp]
    \centering
    \includegraphics[width=0.98\textwidth]{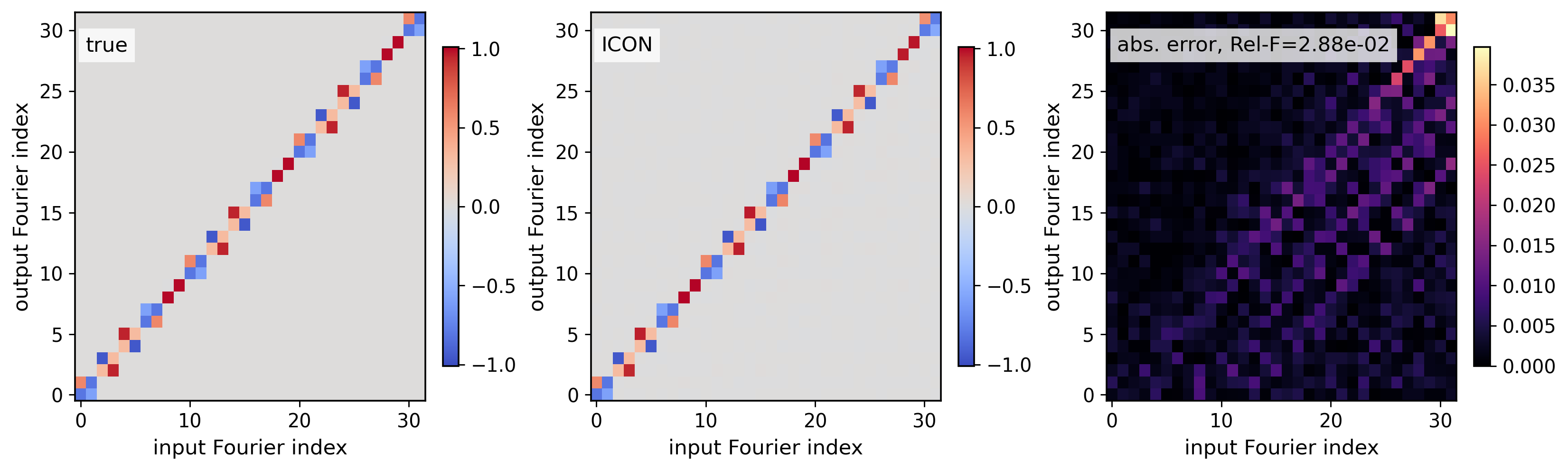}
    \caption{Fourier feature maps for advection at \(c=0.80\) (shift \(s=0.40\)), using the first \(16\) sine--cosine Fourier pairs: exact block-rotation map, ICON Jacobian map, and normalized entrywise error map.}
    \label{fig:adv_feature_map}
\end{figure}

The Fourier feature map in Figure~\ref{fig:adv_feature_map} gives the same
conclusion in matrix form.  The exact map is block diagonal: each frequency has
one sine--cosine rotation block, and there should be no coupling between
different frequencies.  The learned ICON map reproducing this structure is a
stronger statement than low prediction error.  It says that, after reading the
context, the model has induced a local operator with the same algebraic Fourier
structure as the true shift.  Any visible off-block mass should be interpreted
as spurious mode coupling, because the physical advection operator contains no
such mechanism.

The conclusion is that, for the simplest linear operator family, the audit recovers
more than output accuracy.  It verifies amplitude preservation, phase rotation,
and the absence of cross-frequency mixing.  These three properties together are
what characterize the local mechanism of transport.

\FloatBarrier

% ------------------------------------------------------------
\section{Example II: Burgers Equation with Hidden Viscosity}

Burgers dynamics add the first nonlinear mechanism.  We consider
\[
    u_t+u u_x=c u_{xx},\qquad x\in[0,1],
\]
with periodic boundary conditions and hidden viscosity \(c\).  The viscosity
controls how strongly high frequencies are damped, while the nonlinear
convection term makes the tangent map depend on the trajectory around which we
linearize.  This example is therefore designed to separate two effects that a
single prediction error cannot distinguish: context-dependent diffusion and
base-state-dependent mode coupling.

\begin{figure}[!htbp]
    \centering
    \includegraphics[width=0.55\textwidth]{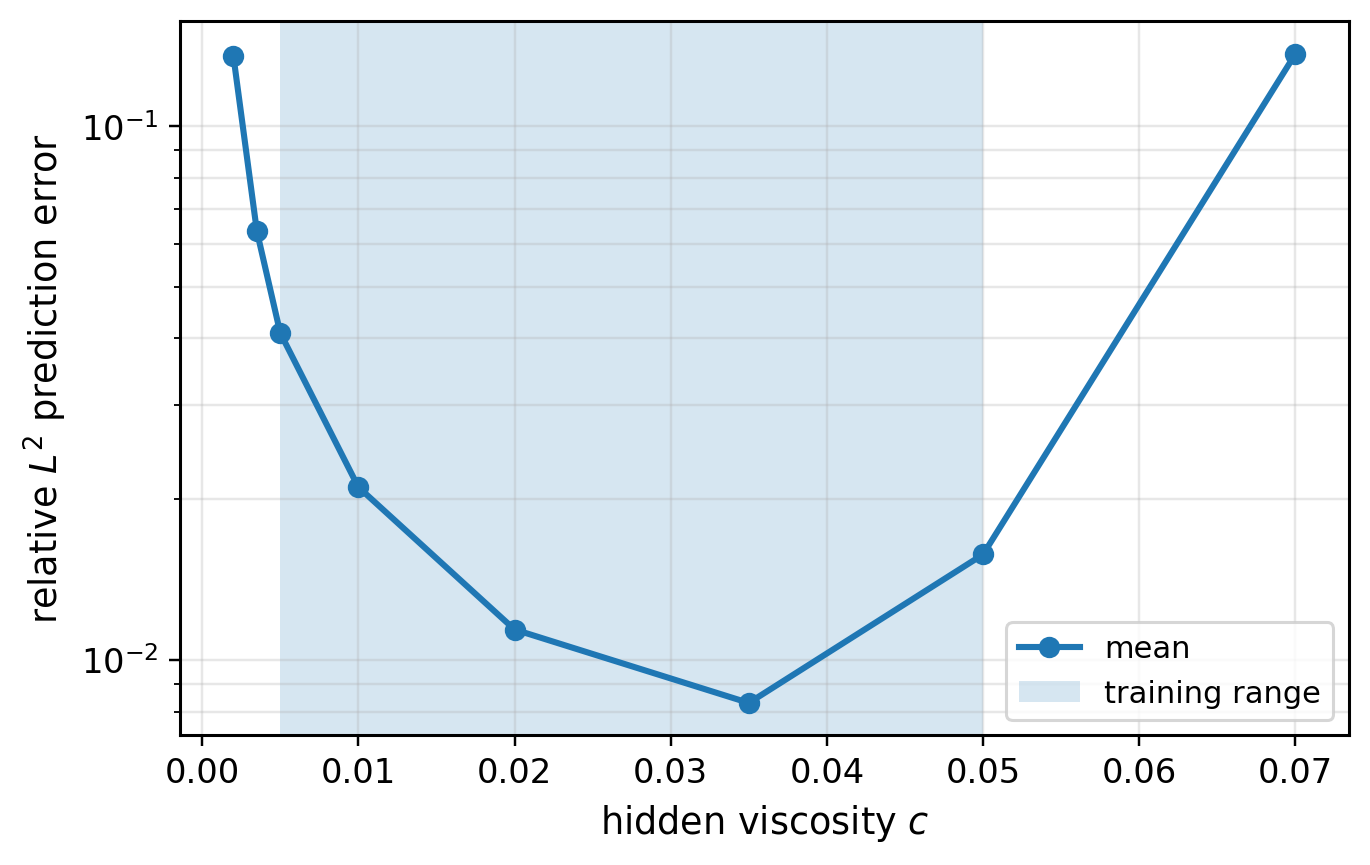}
    \caption{Burgers relative \(L^2\) prediction error with respect to the true finite-time solution \(G_c^T u_0\). The shaded region denotes the training range.}
    \label{fig:burgers_pred_error}
\end{figure}

The prediction curve in Figure~\ref{fig:burgers_pred_error} first compares the
ICON prediction with the true finite-time Burgers solution.  The
relative \(L^2\) error is small on the trained viscosity range and increases
outside it, as expected for very small or large viscosities.  Inside the training
interval, the spectral diagnostics below can therefore be interpreted as errors
in the learned tangent operator, rather than simply as failures of pointwise
solution prediction.

We next audit the context-conditioned query Jacobian at controlled base states.
For a query trajectory \(u(t,x)\), the true tangent perturbation satisfies
\[
    v_t+u v_x+v u_x=c v_{xx},\qquad v(0,x)=v_0(x).
\]
The context controls the hidden viscosity \(c\), which determines the diffusive
term \(c v_{xx}\).  The query base state controls the nonlinear linearization
point through the coefficients \(u\) and \(u_x\) in \(u v_x+v u_x\).  Thus, to
understand the learned Burgers tangent map, we vary the context viscosity while
also checking whether the result changes with the base state.

\begin{figure}[!htbp]
    \centering
    \includegraphics[width=0.8\textwidth]{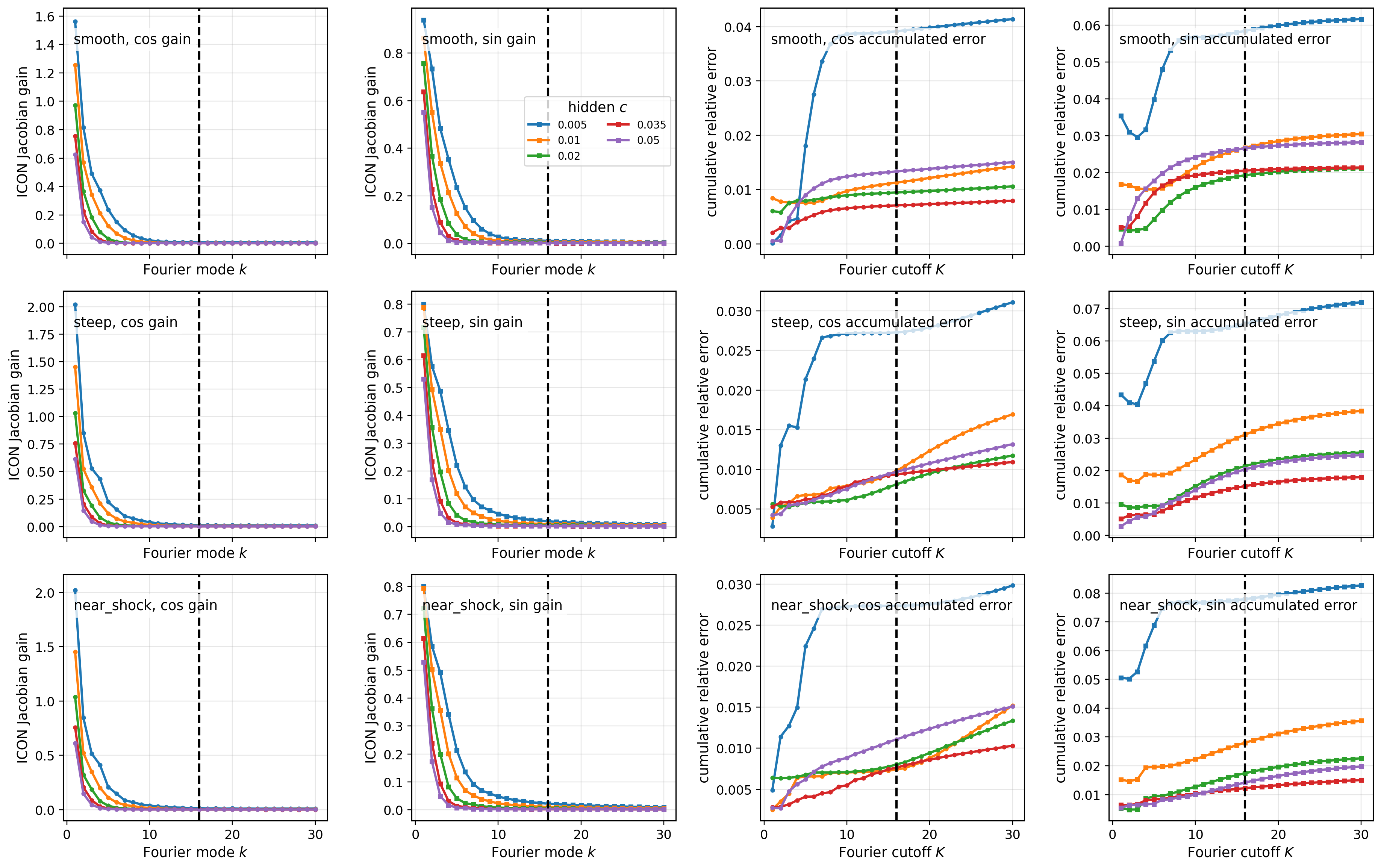}
    \caption{Burgers context-conditioned query Jacobian spectra at deterministic base states. Each row fixes one query base state, and different curves correspond to different context viscosities. The Fourier mode \(k\) denotes the perturbation direction applied to the query input. Gain curves show context-conditioned damping, while error panels report accumulated relative spectral errors with respect to the true Burgers tangent spectrum.}
    \label{fig:burgers_base_context_spectrum}
\end{figure}

In Figure~\ref{fig:burgers_base_context_spectrum}, each row fixes one
deterministic query base state: a smooth profile, a steep-gradient profile, or a
near-shock-like profile.  Within a row, different curves correspond to different
context viscosities.  The Fourier mode \(k\) refers to perturbations of the
query input, not to perturbations of the context functions.  The gain panels show
the learned Jacobian response, while the error panels report the accumulated
relative error against the true Burgers tangent spectrum computed from the
reference solver.  Hence the figure measures not only whether ICON produces the
right qualitative damping order, but also how close its local spectrum is to the
true linearized PDE response.

The viscosity-dependent damping appears as an ordered decay of the gain curves:
smaller viscosity gives weaker damping and larger high-frequency response,
whereas larger viscosity suppresses high modes more strongly.  Differences
between rows reveal the base-state dependence of the nonlinear tangent operator.
The accumulated-error panels show that sharper base states are more demanding,
because stronger gradients in \(u\) and \(u_x\) lead to stronger Fourier-mode
coupling.  Prediction checks for these deterministic probes are reported in
Appendix~\ref{app:base_prediction_checks}.

\begin{figure}[htbp]
    \centering
    \includegraphics[width=0.7\textwidth]{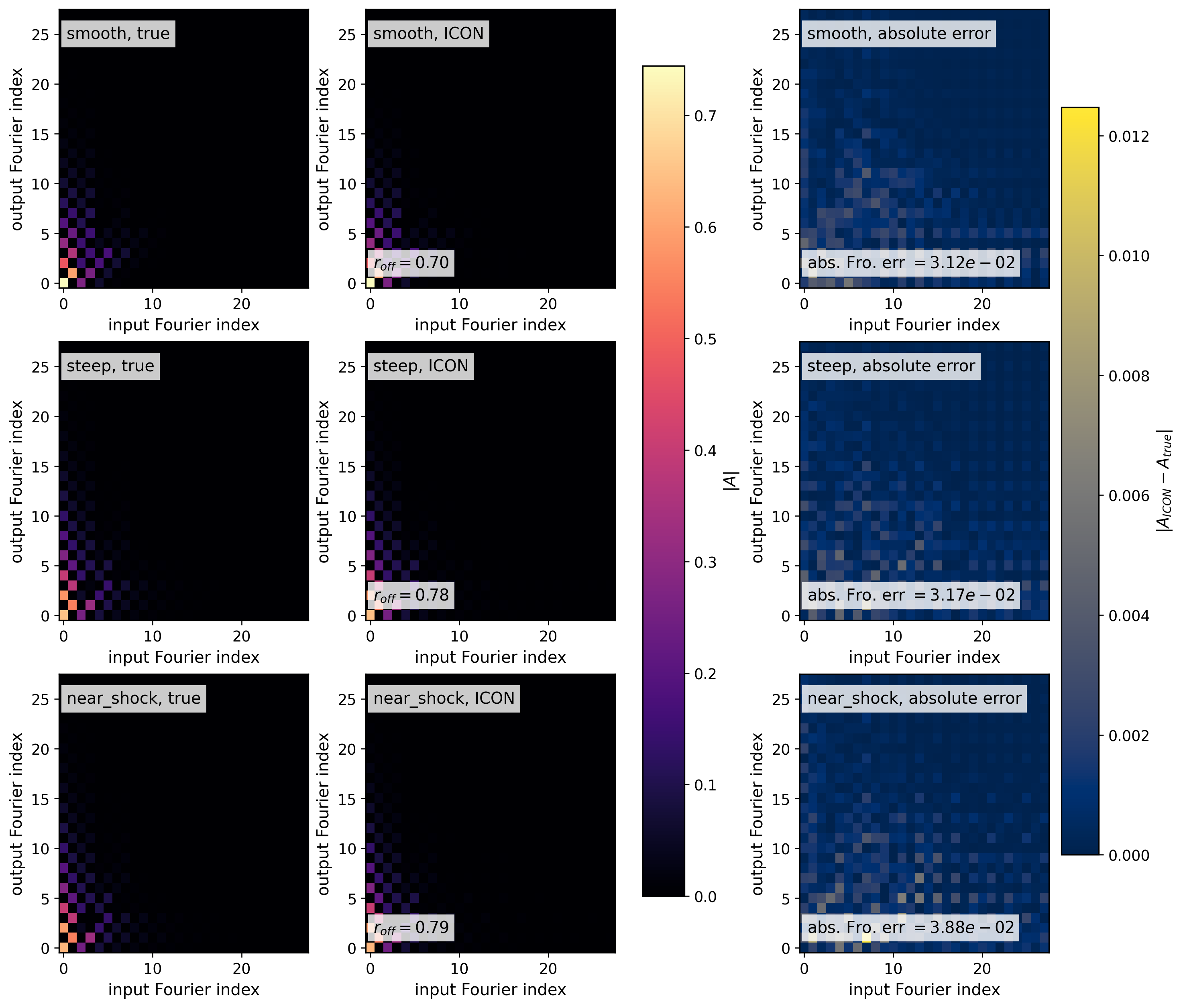}
    \caption{Fourier feature maps for Burgers at \(c=0.02\), evaluated around the smooth, steep-gradient, and near-shock-like base states and projected onto the first \(14\) sine--cosine Fourier pairs. True and ICON maps use a common color scale; error maps show normalized entrywise discrepancy between the learned and true Fourier-projected Jacobians.}
    \label{fig:burgers_feature_map}
\end{figure}

The Fourier feature map in Figure~\ref{fig:burgers_feature_map} gives the
matrix-level view of the same base-state effect.  If the tangent operator were
diagonal in Fourier space, each input mode would mainly excite the same output
mode.  Burgers is different: multiplication by \(u\) and \(u_x\) in
\(u v_x+v u_x\) corresponds to convolution in Fourier space and creates
off-diagonal mode coupling.  The true feature map is obtained by projecting the
reference Burgers tangent operator onto the same Fourier basis as the ICON
Jacobian.  The error panels therefore show the entrywise discrepancy between the
learned Fourier-projected Jacobian and the true one.

Sharper base states have stronger spatial variation, so they should produce
more visible off-diagonal mass.  The feature maps test whether ICON learns this
nonlinear coupling, rather than only a viscosity-dependent smoother.  Agreement
with the true map indicates that the learned local operator captures both the
diagonal damping induced by viscosity and the off-diagonal coupling induced by
the nonlinear base trajectory.

The clean Burgers audit therefore separates two mechanisms. The context viscosity is reflected mainly in the decay of the query-Jacobian gain curves, while the nonlinear base-state effect appears most clearly through off-diagonal Fourier coupling. The accumulated spectral errors and feature-map discrepancies quantify how closely these learned mechanisms match the true Burgers tangent operator. We next use the same Burgers setting as a negative control to test whether the audit also detects a failure of the prompt structure itself.

 In the clean setting, all context solution pairs and the query pair are generated by the same hidden viscosity:
\begin{equation}
\cC_{c}
=
\{u_i,G_c^T u_i\}_{i=1}^{5},
\qquad
y_q=G_c^T u_q .
\end{equation}
In the corrupted setting, the demonstration pairs and the query pair are still valid Burgers input--output examples, but they are generated with different viscosities:
\begin{equation}
\cC_{c_{\mathrm{demo}}}
=
\{u_i,G_{c_{\mathrm{demo}}}^T u_i\}_{i=1}^{5},
\qquad
y_q=G_{c_{\mathrm{query}}}^T u_q,
\qquad
c_{\mathrm{query}}\ne c_{\mathrm{demo}} .
\end{equation}
Thus, the prompt is internally inconsistent: it no longer corresponds to a single hidden Burgers operator. We train ICON on these corrupted prompts and then evaluate the resulting model on clean Burgers prompts.

\begin{figure}[htbp]
    \centering
    \includegraphics[width=0.5\textwidth]{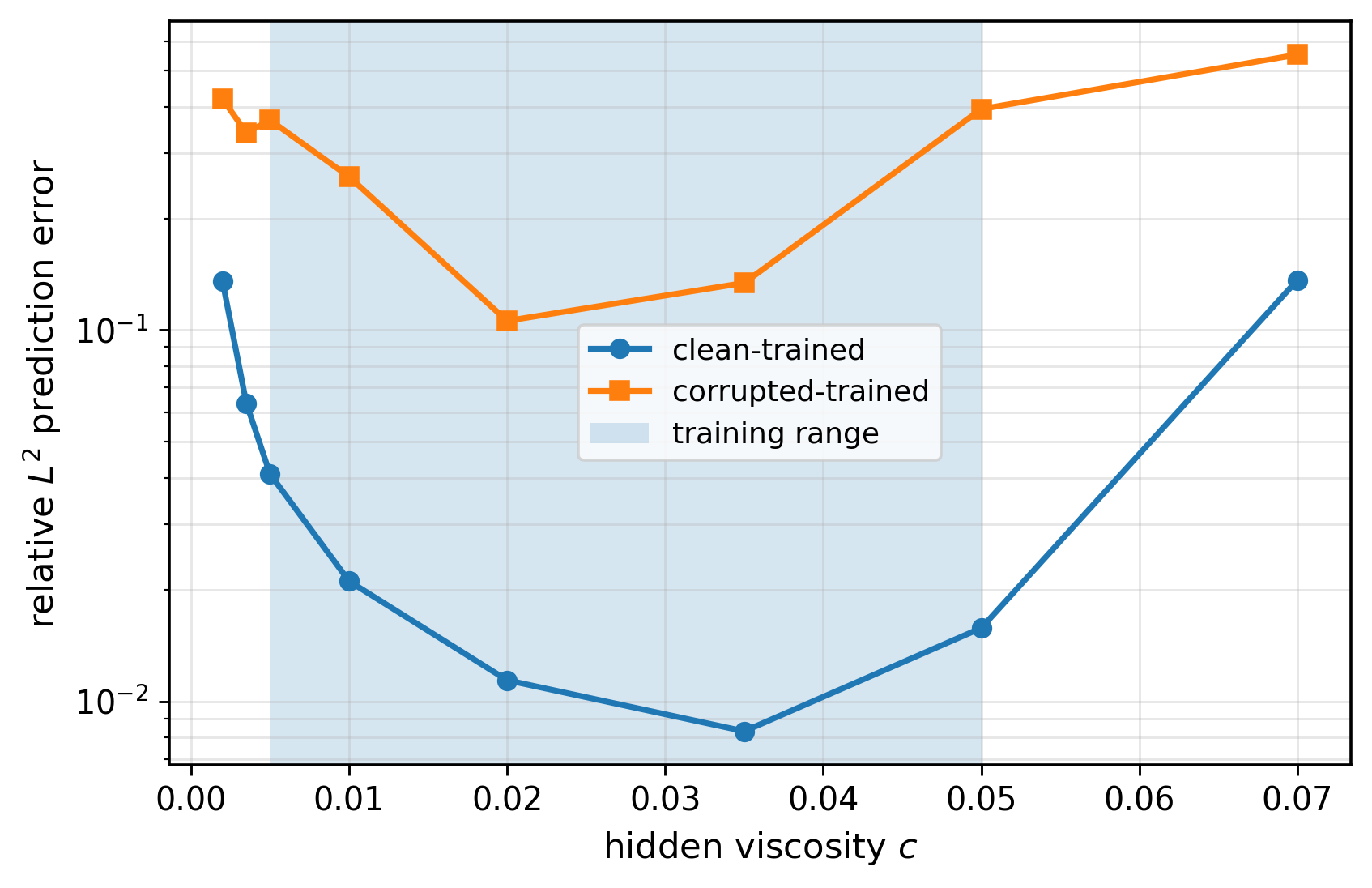}
    \caption{Clean Burgers evaluation of the corrupted-prompt model.  The model
    is trained with inconsistent context--query viscosities, but the test prompts
    use a single correct viscosity \(c\).}
    \label{fig:burgers_wrong_prediction}
\end{figure}

The prediction error in Figure~\ref{fig:burgers_wrong_prediction} shows the most
visible effect of the corrupted training distribution.  But the purpose of this
experiment is not simply to produce a worse predictor.  The more important
question is whether the learned fixed-context map still has the tangent
structure of a clean Burgers operator.  If the context does not identify one
viscosity, the Jacobian should have difficulty matching the true
viscosity-dependent damping profile.

\begin{figure}[htbp]
    \centering
    \includegraphics[width=0.65\textwidth]{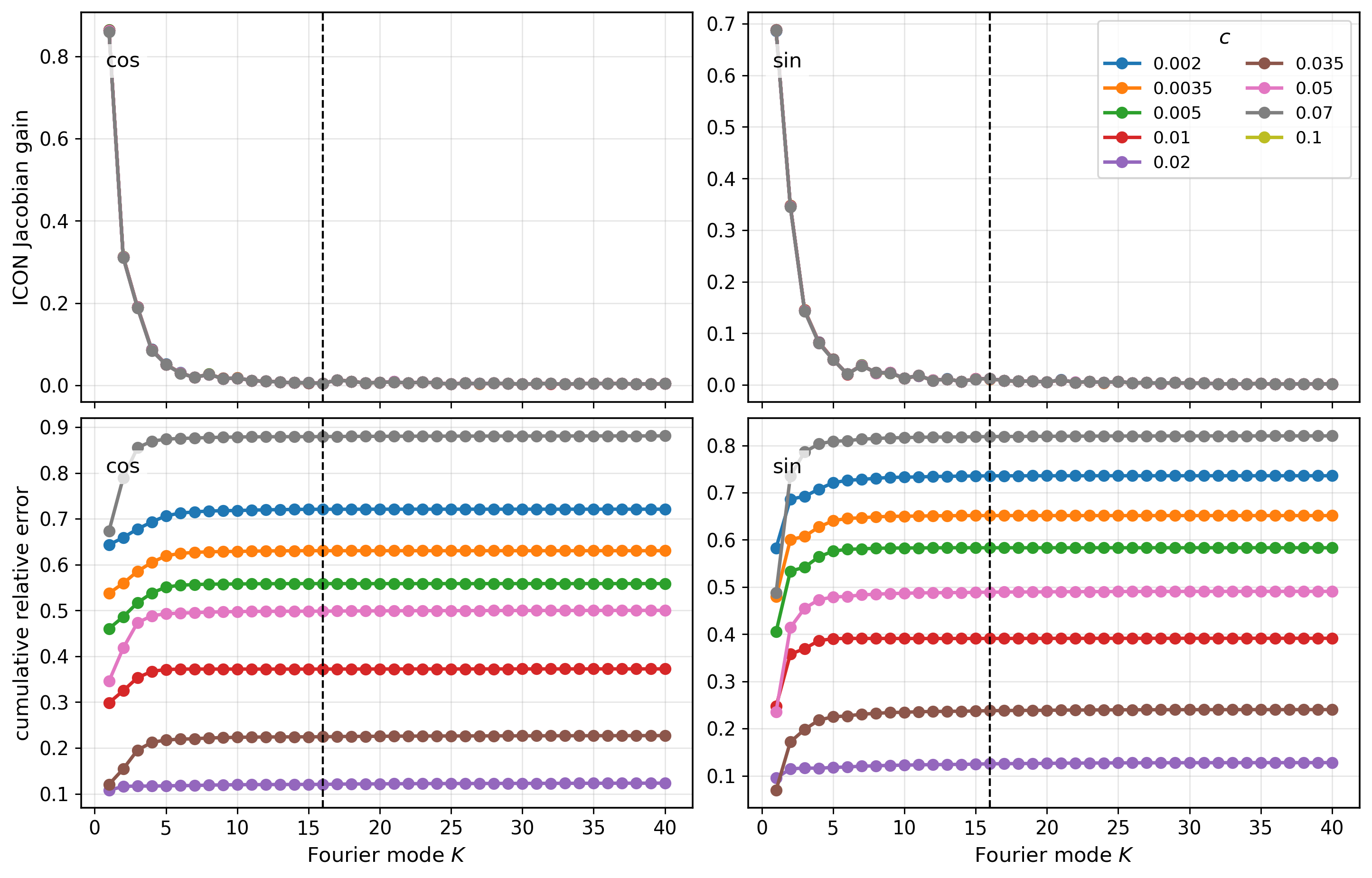}
    \caption{Query-spectrum diagnostic for the corrupted-prompt Burgers model
    evaluated on clean prompts.  The true reference is the clean Burgers tangent
    spectrum, and the error panels report accumulated relative spectral error.}
    \label{fig:burgers_wrong_query_spectrum}
\end{figure}

The spectrum in Figure~\ref{fig:burgers_wrong_query_spectrum} is therefore a
negative-control version of the main Burgers spectrum.  A cleanly trained model
should show an ordered damping profile as viscosity changes.  The corrupted
model has been trained on prompts where the context and query disagree about
which viscosity is active, so a larger accumulated spectral error is evidence of
operator inconsistency.  This is stronger than saying that the output is less
accurate: it says that the prompt no longer induces the correct local tangent
operator.

\begin{figure}[htbp]
    \centering
    \includegraphics[width=0.7\textwidth]{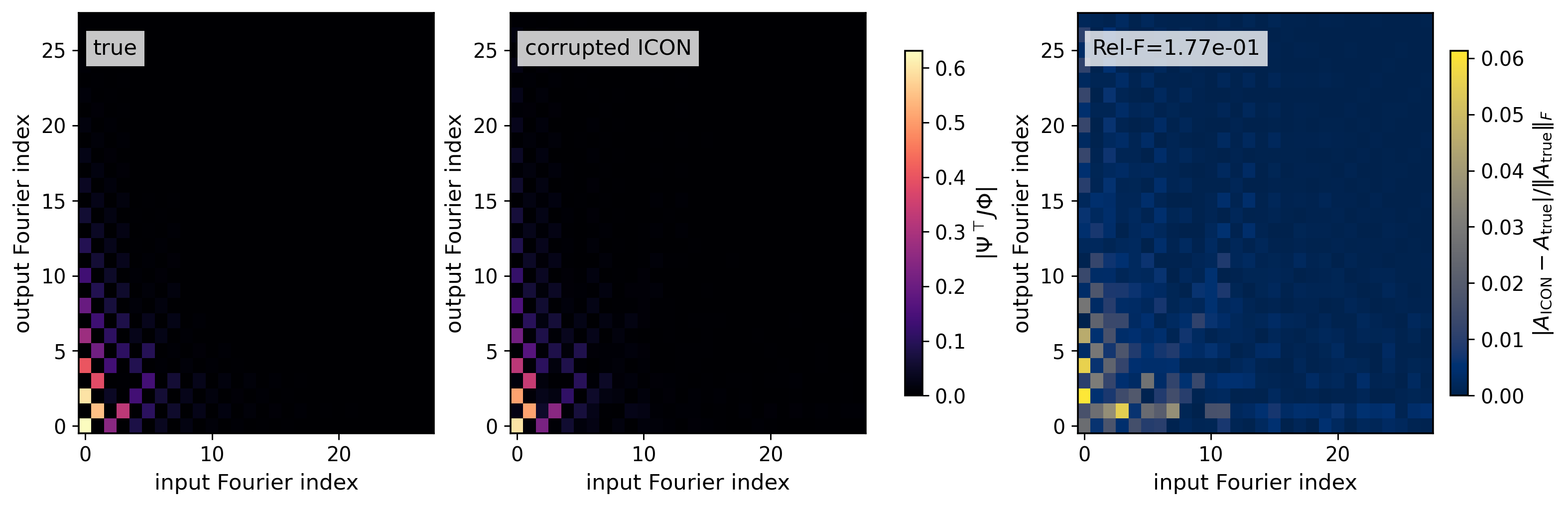}
    \caption{Fourier feature map for the corrupted-prompt Burgers model.  The
    true map is computed from the clean Burgers tangent equation, while the ICON
    map is obtained from the model trained on inconsistent prompts.  The error
    map visualizes the resulting operator-structure mismatch.}
    \label{fig:burgers_wrong_feature_map}
\end{figure}

The feature map in Figure~\ref{fig:burgers_wrong_feature_map} makes the same
point at the level of mode coupling.  In the clean equation, the off-diagonal
structure has a specific source: the linearized convection terms \(u v_x+v u_x\)
around a well-defined trajectory and viscosity.  When the prompt no longer
corresponds to such an operator, the learned Fourier map is not expected to
align with that reference.  This is the intended role of the negative control:
it shows that the spectral audit is sensitive to operator-level consistency in
the prompt, not only to prediction error.

% ------------------------------------------------------------
\section{Example III: Allen--Cahn Equation with Hidden Interface Width}

The Allen--Cahn example tests a different kind of tangent behavior.  We solve
the classical Allen--Cahn equation~\cite{allen1979microscopic}
\[
    u_t=\varepsilon^2u_{xx}+u-u^3,\qquad x\in[0,1],
\]
with hidden parameter \(\varepsilon\).  The learned operator should not behave
like pure diffusion.  The diffusion term damps high frequencies, while the
reaction derivative can amplify low-frequency perturbations near the unstable
state \(u=0\) and suppress perturbations near the stable phases
\(|u|\approx1\).  This makes Allen--Cahn a useful test of whether ICON learns
local stability, not just smoothing.

\begin{figure}[!htbp]
    \centering
    \includegraphics[width=0.5\textwidth]{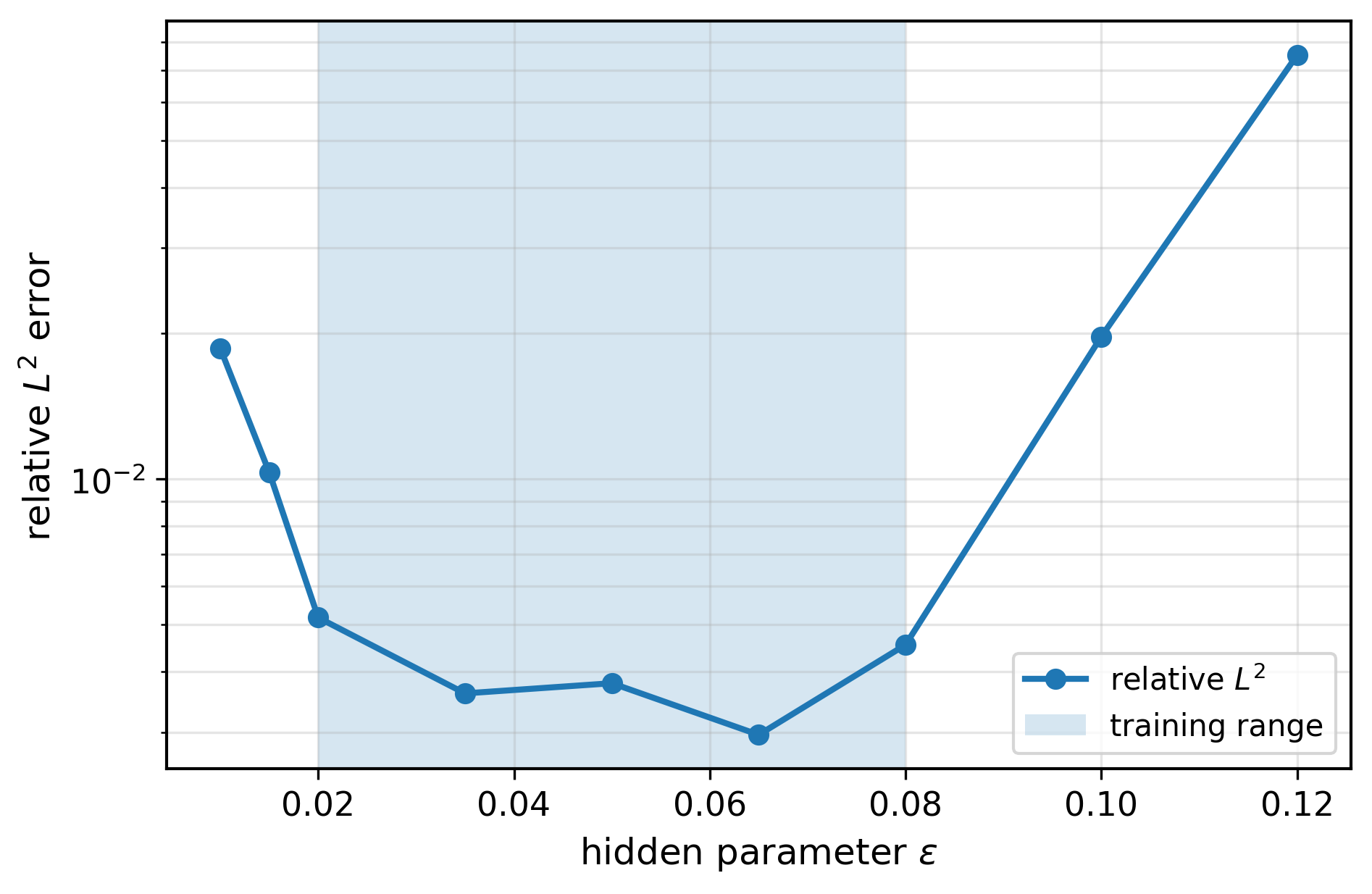}
    \caption{Allen--Cahn relative \(L^2\) prediction error with respect to the true finite-time solution. The shaded region denotes the training range.}
    \label{fig:ac_pred_error}
\end{figure}

The prediction-error curve in Figure~\ref{fig:ac_pred_error} first compares the
ICON prediction with the true finite-time Allen--Cahn solution
\(G_\varepsilon^T u_0\).  The relative \(L^2\) error is small on the training
interval for \(\varepsilon\) and increases outside it.  Thus, inside the training
range, the spectral diagnostics below can be interpreted as errors in the
learned tangent operator, rather than only as failures of pointwise solution
prediction.

We next audit the context-conditioned query Jacobian at controlled base states.
For a trajectory \(u(t,x)\), the true tangent perturbation satisfies
\[
    v_t=\varepsilon^2v_{xx}+(1-3u(t,x)^2)v,
    \qquad v(0,x)=v_0(x).
\]
The role of the hidden parameter is seen from the diffusion part of this tangent equation: on the \(k\)-th Fourier mode, \(v_{xx}\) contributes the multiplier \(-4\pi^2 k^2\), and hence the diffusion term contributes \(-4\pi^2\varepsilon^2 k^2\). Therefore larger \(\varepsilon\) or larger \(k\) leads to stronger damping of high-frequency perturbations, while the reaction coefficient \(1-3u(t,x)^2\) accounts for the base-state-dependent growth and mode coupling.  The query base state controls
the reaction coefficient \(1-3u(t,x)^2\), which determines whether low-frequency
perturbations are locally amplified or damped.  Thus the Allen--Cahn tangent map
depends both on the context-inferred diffusion scale and on the base state around
which the equation is linearized.

\begin{figure}[!htbp]
    \centering
    \includegraphics[width=\textwidth]{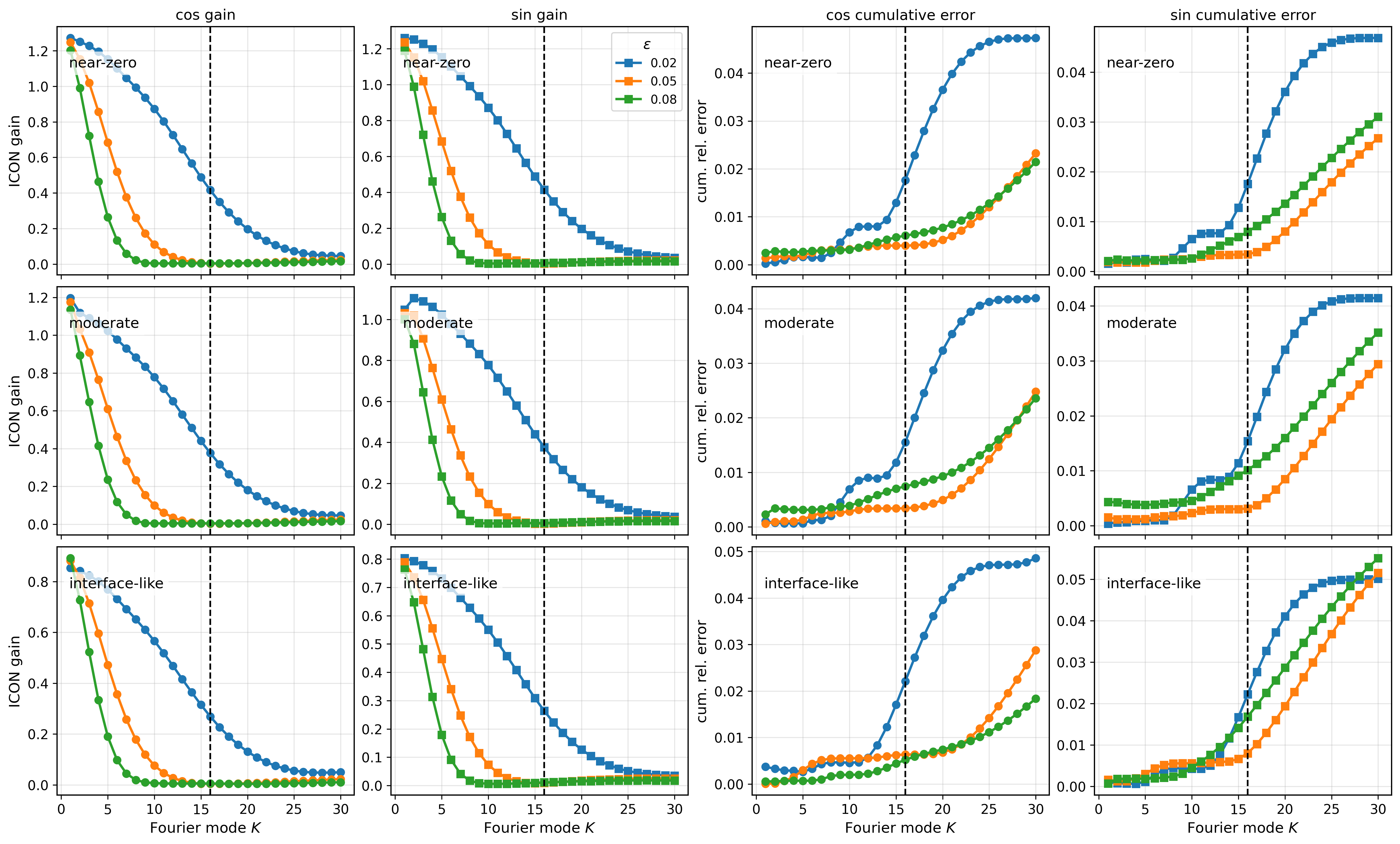}
    \caption{Allen--Cahn context-conditioned query Jacobian spectra at deterministic base states. Each row fixes one query base state, and different curves correspond to different context values of \(\varepsilon\). The Fourier mode \(k\) denotes the perturbation direction applied to the query input. Gain curves show the learned reaction--diffusion response, while error panels report accumulated relative spectral errors with respect to the true Allen--Cahn tangent spectrum.}
    \label{fig:ac_base_context_spectrum}
\end{figure}

With this tangent interpretation in mind, Figure~\ref{fig:ac_base_context_spectrum} compares the learned query Jacobian with the reference Allen--Cahn tangent map. Each row fixes one deterministic query base state: a near-zero sinusoid, a moderate-amplitude sinusoid, or an interface-like profile. Within each row, different curves correspond to different context values of ($\varepsilon$). The Fourier mode (k) refers to perturbations of the query input, not to perturbations of the context functions. The gain panels show the learned Jacobian response, while the error panels report the accumulated relative error against the true tangent spectrum computed from the reference solver.

The main purpose of this test is not to re-identify the Allen--Cahn mechanism, but to check whether the context-conditioned model has learned it at the level of local sensitivities. The reference spectrum is governed by the competition between reaction and diffusion: near the unstable state, broad perturbations may grow, while near the stable phases the response becomes more contractive; increasing ($\varepsilon$) strengthens the damping of high-frequency modes. A correct in-context operator should therefore change its query Jacobian consistently with both the inferred ($\varepsilon$) and the chosen base state. Figure~\ref{fig:ac_base_context_spectrum} shows that this is largely the case: the learned spectra follow the expected shift from low-frequency growth to high-frequency damping, and the accumulated-error panels quantify the remaining discrepancy from the reference tangent solver. Thus, the audit provides information that is not contained in the prediction error alone: it tests whether ICON has learned a locally meaningful Allen--Cahn operator, rather than only matching deterministic query outputs.

\begin{figure}[H]
    \centering
    \includegraphics[width=0.6\textwidth]{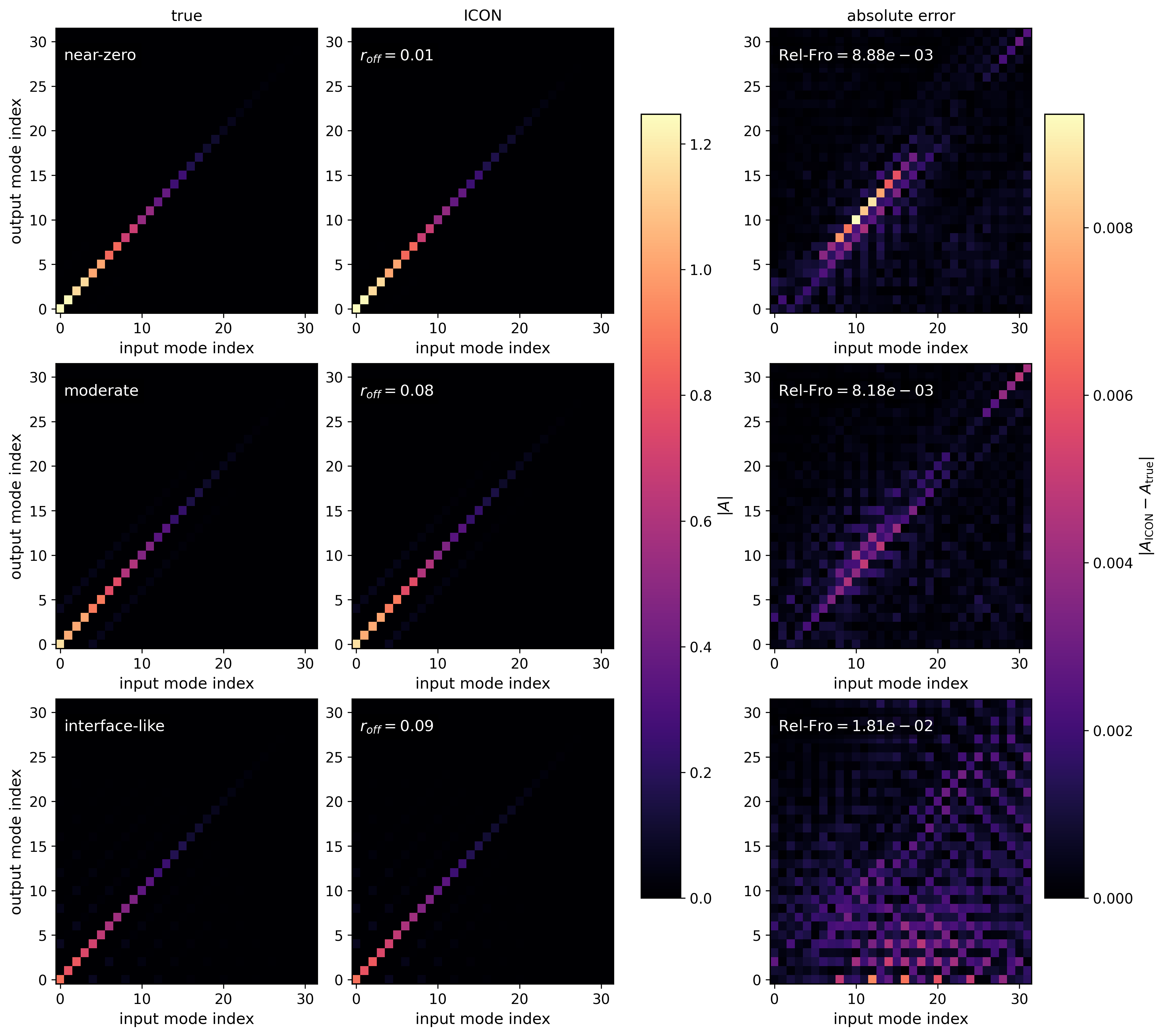}
    \caption{Fourier feature maps for Allen--Cahn at \(\varepsilon=0.05\), evaluated around the near-zero, moderate-amplitude, and interface-like base states and projected onto the first \(16\) sine--cosine Fourier pairs. True and ICON maps use a common color scale; error maps show normalized entrywise discrepancy between the learned and true Fourier-projected Jacobians.}
    \label{fig:ac_feature_map}
\end{figure}

Figure~\ref{fig:ac_feature_map} gives a matrix-level check of the learned Allen--Cahn tangent operator. The main observation is that both the reference map and the ICON map are largely diagonal-dominant, which is consistent with the fact that diffusion acts diagonally in Fourier space and is the dominant source of high-frequency damping. The remaining off-diagonal entries come from the reaction coefficient $(1-3u(x)^2)$, whose spatial variation couples nearby modes. Thus, unlike Burgers, the Allen--Cahn tangent map should not exhibit strong global mode mixing.

The feature map shows that ICON captures this qualitative structure: it preserves the dominant diagonal response and does not introduce large spurious off-diagonal coupling. The error map is therefore useful not simply as another accuracy plot, but as a check of whether the learned local operator has the correct matrix structure. In this example, agreement with the reference map indicates that the model has learned the main stability mechanism of Allen--Cahn: high-frequency damping controlled by diffusion, together with weaker base-state-dependent coupling from the reaction term.

The Allen--Cahn audit therefore tests a local stability mechanism rather than
only smoothing.  The context value of \(\varepsilon\) is seen mainly through
high-frequency damping, while the base state controls whether low modes are
amplified or suppressed through the reaction derivative.  The accumulated
spectral errors and feature-map discrepancies quantify how far these learned
mechanisms are from the true Allen--Cahn tangent operator.

\FloatBarrier

% ------------------------------------------------------------
\section{Discussion}

The experiments show that the same Jacobian-based audit can reveal different
PDE mechanisms.  In advection, the learned operator should preserve amplitudes,
rotate Fourier phase, and avoid cross-frequency mixing.  In Burgers, the audit
separates viscosity-dependent damping from nonlinear mode coupling.  In
Allen--Cahn, it captures the competition between reaction-driven low-frequency
growth and diffusion-driven high-frequency decay.  These examples show that
prediction accuracy and local operator structure are different aspects of the
learned model.

The diagnostics play complementary roles.  Prediction error checks whether the
output is close to the reference solution.  Gain curves and accumulated spectral
errors compare the learned query Jacobian with the true tangent map.  Phase and
block errors are essential for transport, while Fourier feature maps reveal
mode coupling that scalar gains can hide.  The ablations in
Appendix~\ref{app:ablation_studies} further show that the audit is sensitive to
training bandwidth and to prompt inconsistency.

Several extensions are natural.  Dispersive equations would test
frequency-dependent phase speed, while pattern-forming systems such as
Cahn--Hilliard dynamics~\cite{cahn1958free} would provide a more delicate
mode-stability structure.  In higher dimensions or on irregular geometries, the
Fourier basis could be replaced by localized wave packets, Laplacian
eigenfunctions, graph bases, or randomized tangent directions.  The same audit
could also be used to compare ICON with other conditional neural operators, such
as parameter-conditioned DeepONet/FNO models, conditional neural
processes~\cite{garnelo2018conditional}, or continuous-depth dynamics
models~\cite{chen2018neural}.

% ------------------------------------------------------------
\section{Conclusion}

We introduced a spectral audit for in-context operator learning.  After
a context is fixed, ICON defines a query-to-output map.  By differentiating this
map with respect to the query and probing the Jacobian in Fourier directions, we
obtain a local view of the operator inferred from the prompt.

Across the three PDE examples, the learned Jacobians reproduce the expected
mechanisms: phase rotation for advection, viscosity-dependent damping and
nonlinear mode coupling for Burgers, and the competition between reaction and
diffusion in Allen--Cahn.  The ablation and negative-control experiments clarify
what these diagnostics are sensitive to: training bandwidth matters for
high-frequency phase, and prompt inconsistency damages the learned tangent
operator.

\section*{Acknowledgments}

This work was supported in part by the DARPA-DIAL grant HR00112490484, the ONR Vannevar Bush Faculty Fellowship under grant N00014-22-1-2795, and, for Liu Yang, the National Research Foundation, Singapore, under the NRF Fellowship, Project No. NRF-NRFF17-2025-0006. The authors would also like to thank Professor Yannis Kevrekidis for useful discussions that helped initiate this research.

\appendix

\section{Auxiliary Prediction Checks for Deterministic Base States}
\label{app:base_prediction_checks}

The main text uses deterministic base states to probe the learned Jacobian. Since these base states are not sampled from the random training distribution, we first check whether ICON predicts their final-time solutions accurately.

\begin{figure}[!htbp]
\centering
\includegraphics[width=0.7\textwidth]{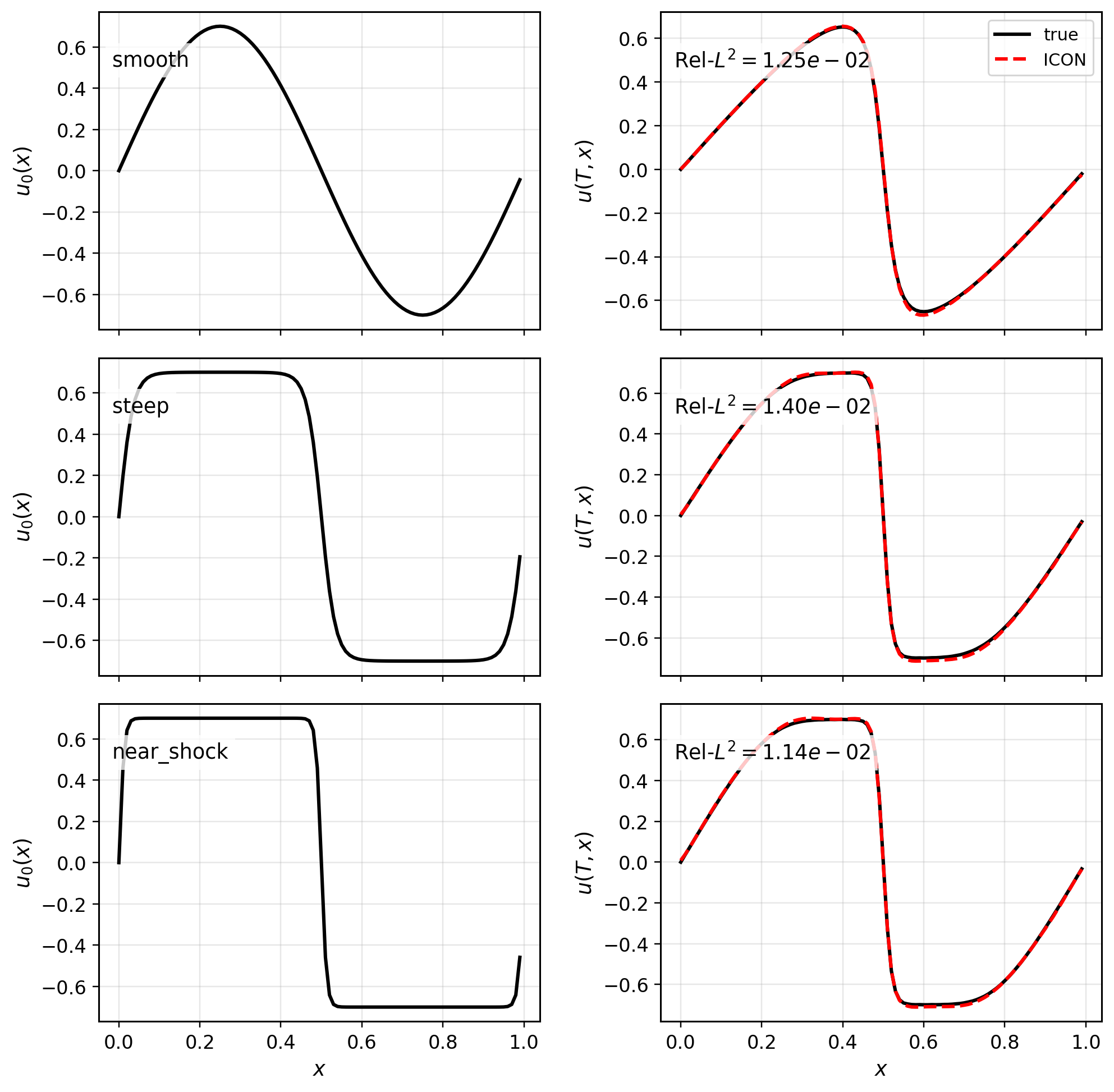}
\caption{Burgers predictions for the smooth, steep-gradient, and near-shock-like base states.}
\label{fig:burgers_base_pred}
\end{figure}

Figure~\ref{fig:burgers_base_pred} shows that ICON gives accurate predictions for the three Burgers base states used in the Jacobian audit.

\begin{figure}[!htbp]
\centering
\includegraphics[width=0.7\textwidth]{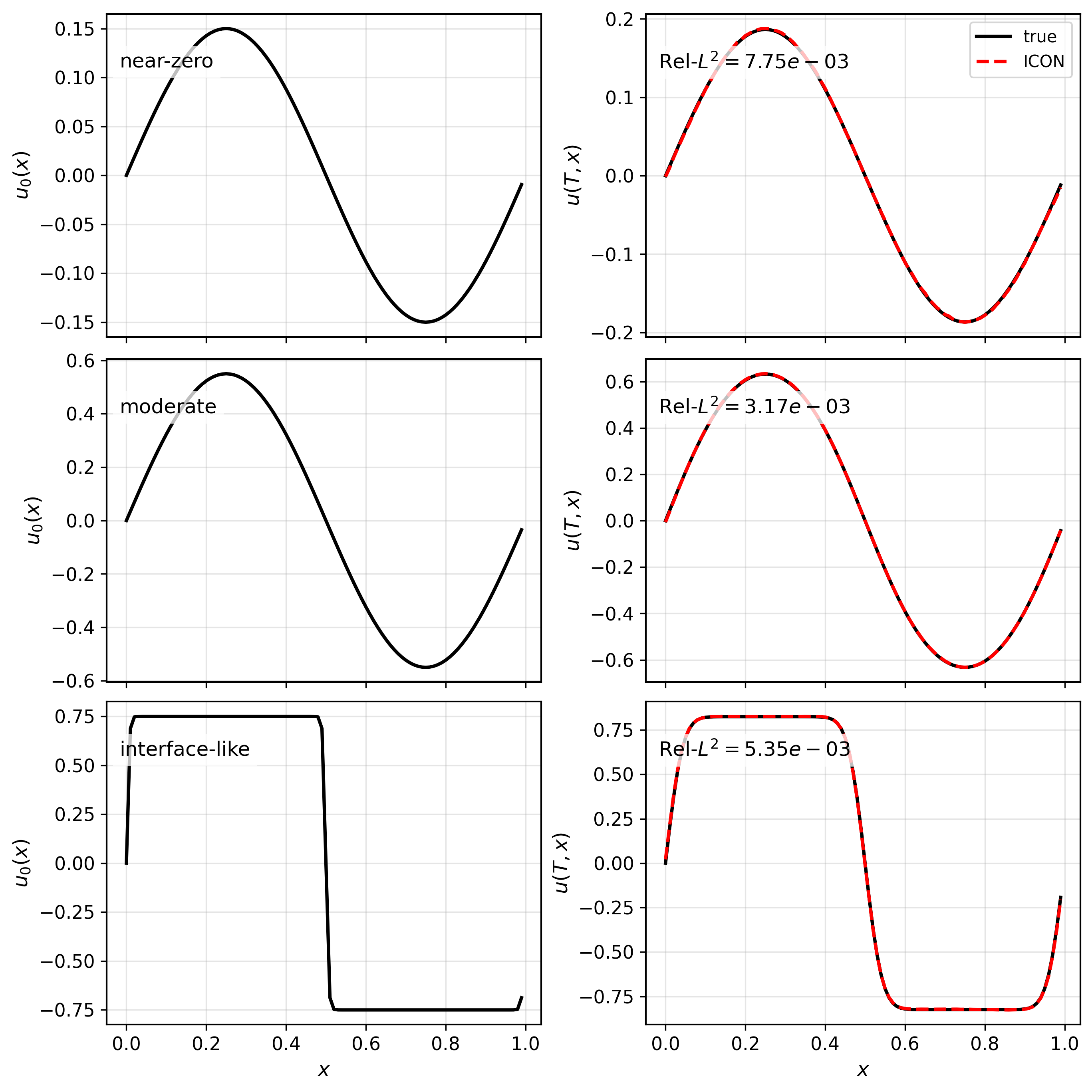}
\caption{Allen--Cahn predictions for the near-zero, moderate-amplitude, and interface-like base states.}
\label{fig:ac_base_pred}
\end{figure}
Figure~\ref{fig:ac_base_pred} shows the corresponding prediction checks for the Allen--Cahn base states. These results confirm that the deterministic probes used in the main text are not poorly predicted out-of-distribution states.

\section{Advection bandwidth: training on 8 versus 16 modes}
\label{app:ablation_studies}

For a shift operator, high frequencies are unforgiving.  The phase angle is
proportional to the Fourier mode, so the same shift error is much more visible
at \(k=16\) than at \(k=4\).  We therefore train one advection model on random
initial functions supported up to eight Fourier modes and another on functions
supported up to sixteen modes. Other details are unchanged. 

\begin{figure}[!htbp]
    \centering
    \includegraphics[width=0.7\textwidth]{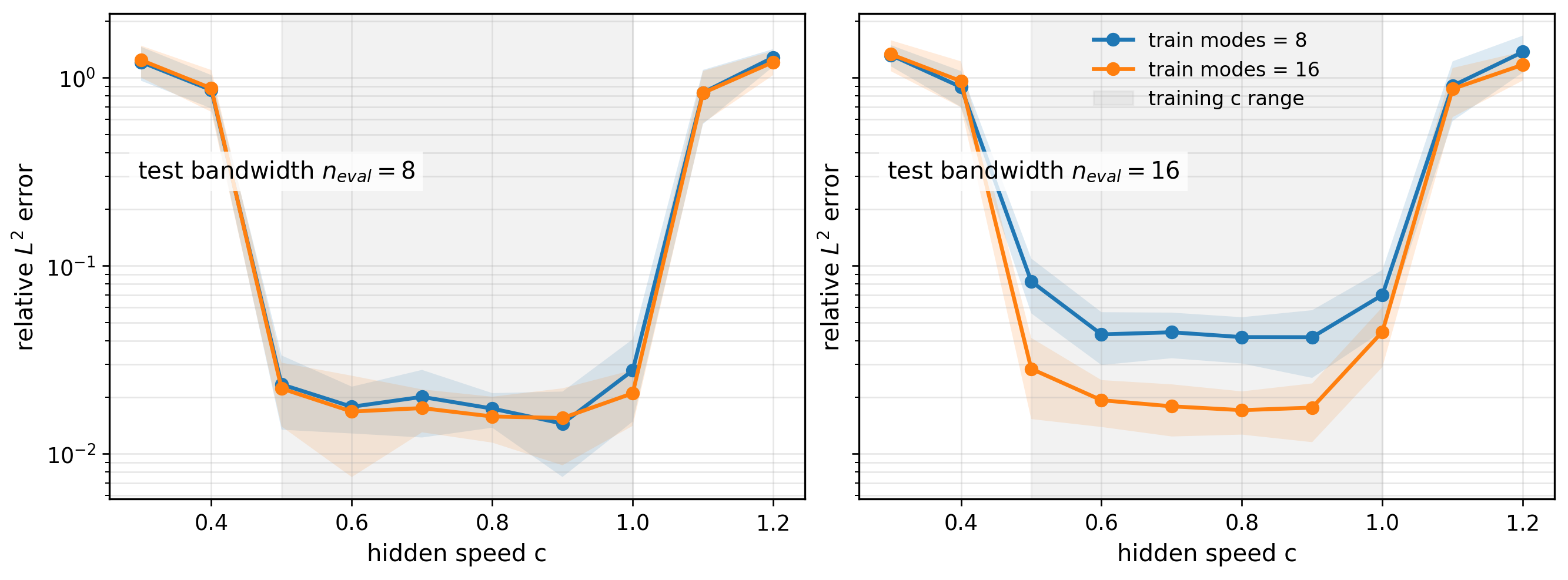}
    \caption{Advection bandwidth ablation.  The two models are trained with
    \(n_{\mathrm{train}}=8\) and \(n_{\mathrm{train}}=16\), respectively, and are
    evaluated on test functions with bandwidth \(n_{\mathrm{eval}}=8\) and
    \(n_{\mathrm{eval}}=16\).  The shaded region denotes the training range for
    the hidden speed \(c\).}
    \label{fig:adv_ablation_prediction}
\end{figure}

Figure~\ref{fig:adv_ablation_prediction} separates interpolation within the training bandwidth from generalization to unseen Fourier modes. When the test functions use ($n_{\mathrm{eval}}=8$), both checkpoints are evaluated on frequencies covered by the 8-mode training data, so their prediction errors are expected to be close. When the test functions use ($n_{\mathrm{eval}}=16$), the 8-mode model must extrapolate the shift rule to higher modes that it has not seen during training, while the 16-mode model has direct training support on those modes. The resulting error gap shows that training bandwidth affects the learned response, and the spectral diagnostics below identify how this degradation appears at the operator level.

\begin{figure}[!htbp]
    \centering
    \includegraphics[width=0.7\textwidth]{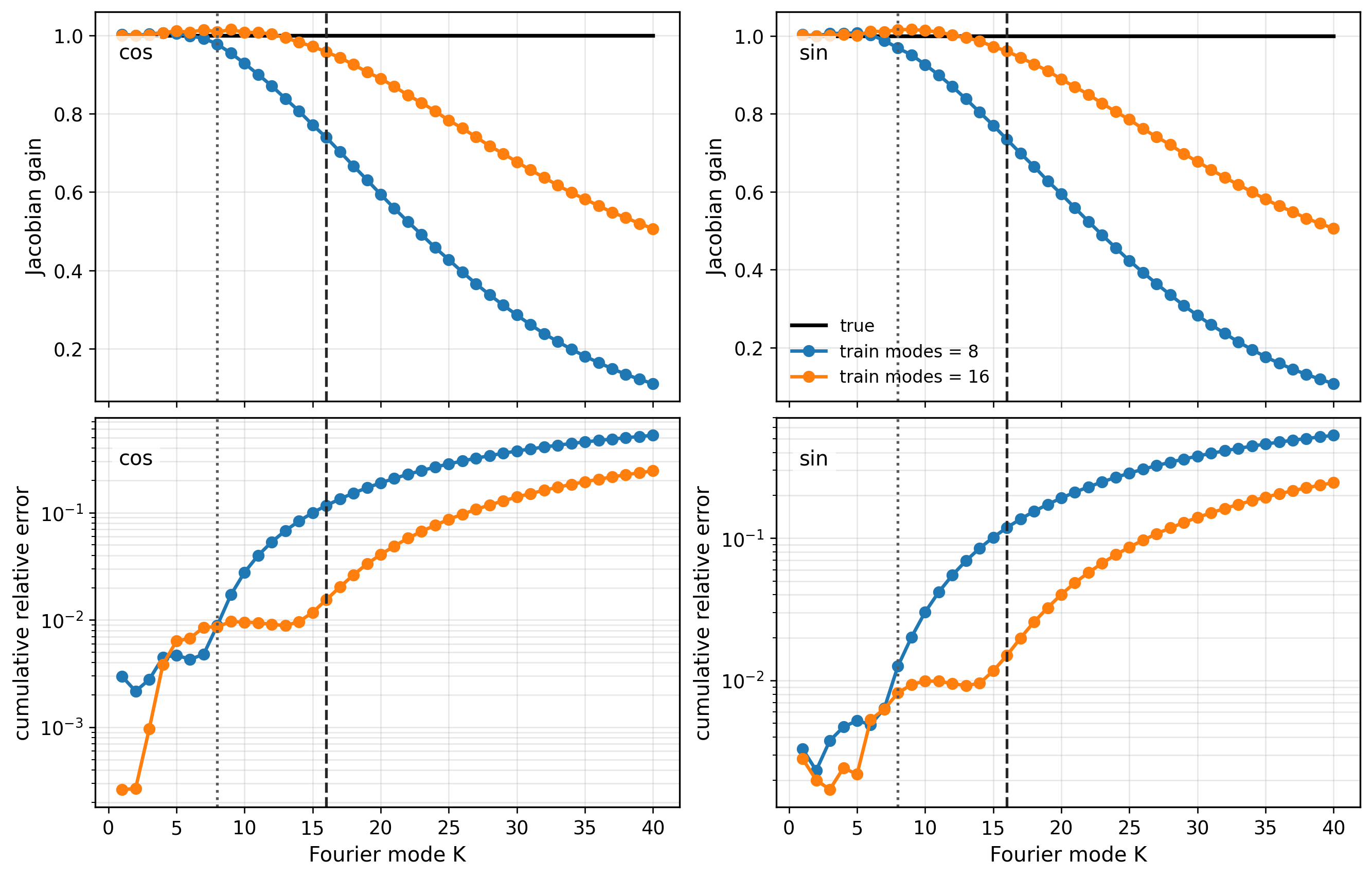}
    \caption{Advection bandwidth ablation at \(c=0.80\).  The plot compares
    the query Fourier spectra and accumulated relative spectral errors of the
    models trained with 8 and 16 Fourier modes.}
    \label{fig:adv_ablation_query_spectrum}
\end{figure}

The query-spectrum ablation in Figure~\ref{fig:adv_ablation_query_spectrum}
turns this into a local statement about the Jacobian.  The true shift has unit
gain at every resolved mode.  Thus, if the 8-mode model begins to deviate beyond
\(k=8\), the issue is not that the shift operator changes there; it is that the
model has not learned to extend the same tangent response to frequencies absent
from its training functions.  This is exactly the sort of distinction that is
hard to make from a single relative \(L^2\) curve.

\begin{figure}[!htbp]
    \centering
    \includegraphics[width=0.5\textwidth]{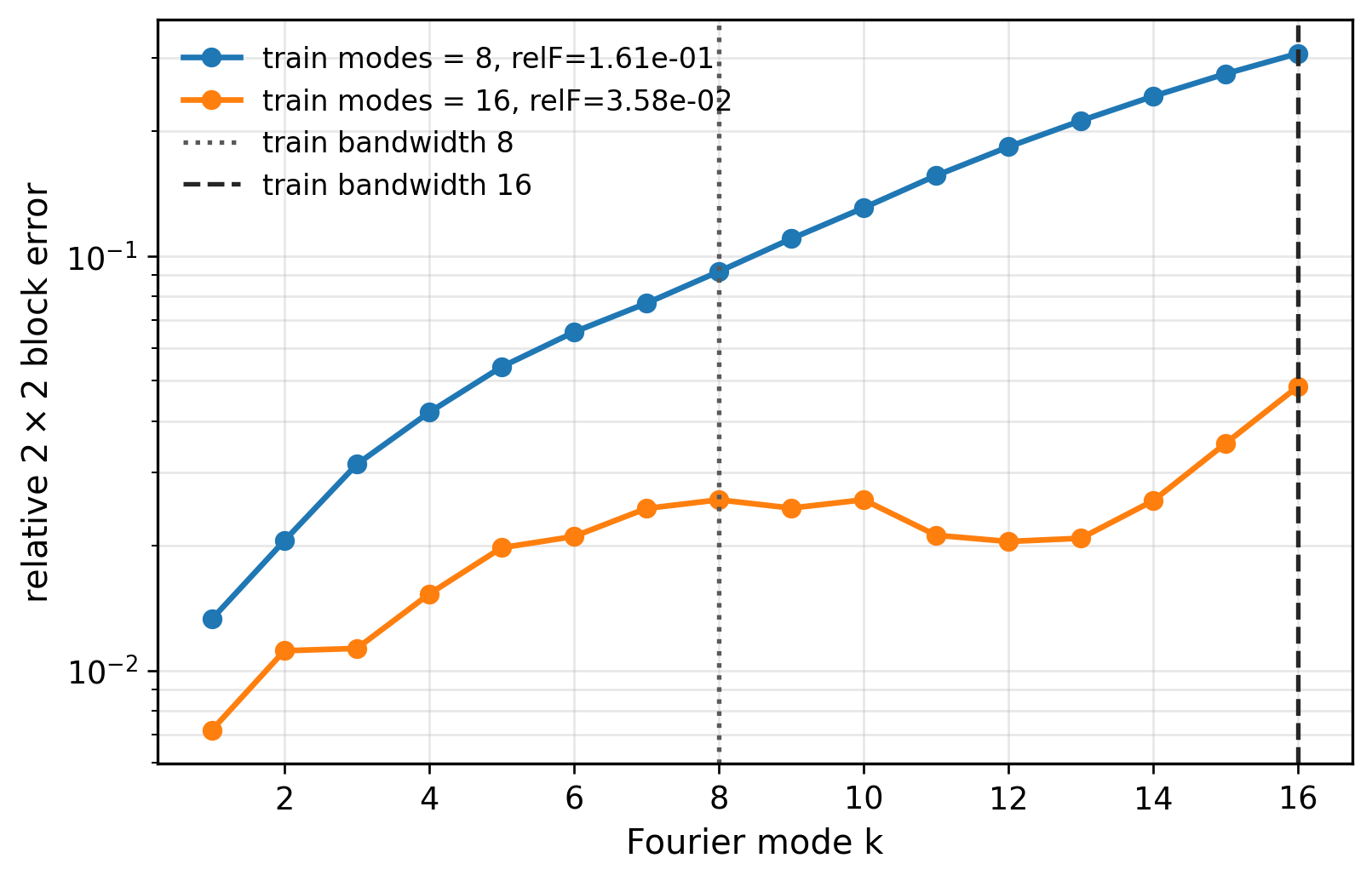}
    \caption{Phase and block-rotation errors for the advection bandwidth
    ablation at \(c=0.80\).  Because the true advection Jacobian is a direct sum
    of \(2\times2\) Fourier rotation blocks, this diagnostic directly measures
    the learned shift phase.}
    \label{fig:adv_ablation_phase}
\end{figure}

\begin{figure}[!htbp]
    \centering
    \includegraphics[width=0.98\textwidth]{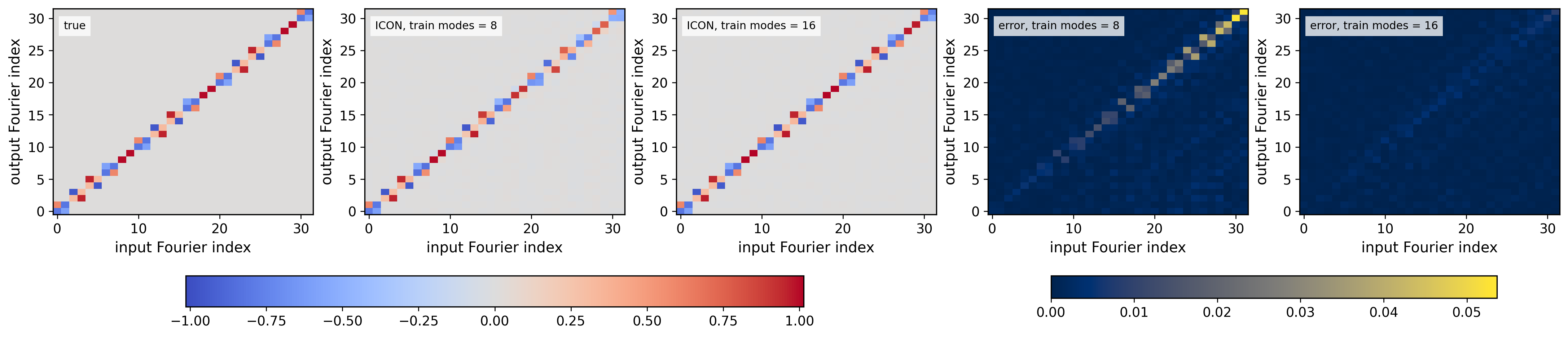}
    \caption{Fourier feature-map ablation for advection at \(c=0.80\).  The
    exact shift map is compared with ICON models trained on 8 and 16 Fourier
    modes.  The projection uses the first \(16\) sine--cosine Fourier pairs, so
    modes beyond \(k=8\) test frequency extrapolation for the 8-mode model.}
    \label{fig:adv_ablation_feature_map}
\end{figure}

Figure~\ref{fig:adv_ablation_phase} shows that the 8-mode model does not only lose accuracy at high frequencies; it loses the correct sine--cosine phase rotation. The 16-mode model preserves the block structure over a wider range, whereas the 8-mode model deteriorates beyond its training bandwidth. The Fourier maps in Figure~\ref{fig:adv_ablation_feature_map} show the same effect visually: low-frequency blocks are learned by both models, but the lower-bandwidth model fails to extrapolate the rotation pattern beyond (k=8). Thus the audit detects a structural generalization failure, not just a larger prediction error.

\FloatBarrier

% ------------------------------------------------------------

\section{True Tangent Maps for the Reference Operators}
\label{app:true_derivatives}

This appendix derives the reference tangent maps used in the spectral audit.
For each PDE family, the true spectrum is computed from the derivative of the
finite-time solution operator with respect to the initial condition \(u_0\).
This matches the learned quantity audited in the main text: the context is held
fixed, and only the query initial condition is perturbed.

Let \(G_\eta^T\) denote the finite-time solution operator at time \(T\), where
\(\eta\) is the hidden parameter.  For an initial perturbation \(v_0\), the true
tangent map is
$$
    D_{u_0}G_\eta^T(u_0)v_0
    =
    \left.\frac{\dd}{\dd\rho}
    G_\eta^T(u_0+\rho v_0)
    \right|_{\rho=0}.
$$
In the implementation, the true spectra are obtained by Jacobian--vector
products through the discrete reference solver.  The equations below give the
continuous variational equations corresponding to these discrete derivatives.

\subsection{Advection}

The periodic advection equation is
$$
    u_t+c u_x=0,
    \qquad u(0,x)=u_0(x).
$$
At time \(T\),
$$
    G_c^{T}u_0(x)=u_0(x-cT),
$$
with the argument understood modulo one.

Let \(u_0^\rho=u_0+\rho v_0\).  Then
$$
    G_c^{T}u_0^\rho(x)
    =
    u_0(x-cT)+\rho v_0(x-cT),
$$
so
$$
    D_{u_0}G_c^{T}(u_0)v_0(x)
    =
    v_0(x-cT).
$$
This derivative is independent of \(u_0\) because the advection operator is
linear.

For the Fourier basis, using
$$
    \cos(2\pi k(x-s))
    =
    \cos(2\pi kx)\cos(2\pi ks)
    +
    \sin(2\pi kx)\sin(2\pi ks),
$$
$$
    \sin(2\pi k(x-s))
    =
    \sin(2\pi kx)\cos(2\pi ks)
    -
    \cos(2\pi kx)\sin(2\pi ks),
$$
with \(s=c T\), the exact real-Fourier block is a rotation.  Thus the
continuous shift has unit gain for every mode, and the phase or block-rotation
error captures the hidden shift.

\subsection{Burgers}

The viscous Burgers equation is
$$
    u_t+u u_x=c u_{xx},
    \qquad u(0,x)=u_0(x).
$$
Let \(u(t,x)\) be the solution generated by \(u_0\).  Perturb the initial
condition by \(u_0^\rho=u_0+\rho v_0\), and let \(u^\rho\) solve
$$
    u_t^\rho+u^\rho u_x^\rho=c u_{xx}^\rho,
    \qquad u^\rho(0,x)=u_0(x)+\rho v_0(x).
$$
Define
$$
    v(t,x)=
    \left.
    \frac{\partial u^\rho(t,x)}{\partial\rho}
    \right|_{\rho=0}.
$$
Differentiating the PDE gives
$$
    v_t+u v_x+v u_x=c v_{xx},
    \qquad v(0,x)=v_0(x).
$$
Therefore
$$
    D_{u_0}G_c^T(u_0)v_0=v(T,\cdot).
$$
The term \(c v_{xx}\) damps high frequencies, while the terms
\(u v_x+v u_x\) create Fourier-mode coupling through multiplication by
spatially varying coefficients.

\subsection{Allen--Cahn}

The Allen--Cahn equation is
$$
    u_t=\varepsilon^2u_{xx}+u-u^3,
    \qquad u(0,x)=u_0(x).
$$
Let \(u(t,x)\) be the solution generated by \(u_0\).  Perturb the initial
condition by \(u_0^\rho=u_0+\rho v_0\), and define
$$
    v(t,x)=
    \left.
    \frac{\partial u^\rho(t,x)}{\partial\rho}
    \right|_{\rho=0}.
$$
Differentiating
$$
    u_t^\rho=\varepsilon^2u_{xx}^\rho+u^\rho-(u^\rho)^3
$$
yields
$$
    v_t=\varepsilon^2v_{xx}+(1-3u(t,x)^2)v,
    \qquad v(0,x)=v_0(x).
$$
Therefore
$$
    D_{u_0}G_\varepsilon^T(u_0)v_0=v(T,\cdot).
$$
If the base trajectory is approximately constant, \(u(t,x)\approx a\), then the
\(k\)-th Fourier mode has approximate growth rate
$$
    \lambda_k(a,\varepsilon)
    =
    1-3a^2-4\pi^2\varepsilon^2k^2.
$$
This explains why low modes can be amplified near \(a=0\), while high modes are
damped by diffusion.

% ------------------------------------------------------------
\bibliographystyle{unsrtnat}
\bibliography{references}

\end{document}